\documentclass[aop]{imsart}

\RequirePackage{amsthm,amsmath,amsfonts,amssymb,mathtools}
\RequirePackage[numbers]{natbib}
\RequirePackage{graphicx}
\RequirePackage[colorlinks,citecolor=blue,urlcolor=blue,linkcolor=blue]{hyperref}

\makeatletter\def\journal@name{}\makeatother
\startlocaldefs
\numberwithin{equation}{section}
\theoremstyle{plain}
\newtheorem{mainthm}{Theorem}
\newtheorem{proposition}{Proposition}[section]
\newtheorem{lemma}[proposition]{Lemma}
\newtheorem{corollary}[proposition]{Corollary}
\theoremstyle{definition}
\newtheorem{definition}[proposition]{Definition}
\newtheorem{remark}[proposition]{Remark}

\newcommand{\R}{\mathbb{R}}
\newcommand{\E}{\mathbb{E}}
\newcommand{\PP}{\mathbb{P}}
\newcommand{\Id}{I}
\newcommand{\dd}{\mathrm{d}}
\newcommand{\Ent}{\mathrm{Ent}}
\newcommand{\Ex}{\mathcal{E}_x}
\newcommand{\Ev}{\mathcal{E}_v}
\newcommand{\Ccor}{\mathcal{C}}
\newcommand{\Leb}{\mathcal{L}}
\newcommand{\Exp}{\mathrm{Exp}}
\newcommand{\Etrap}{E_{\mathrm{trap}}}
\newcommand{\cgap}{\mathrm{gap}}
\newcommand{\Emod}{\mathcal{H}_G}
\newcommand{\Ly}{\mathcal{V}}
\newcommand{\Dst}{\mathcal{D}_*}
\DeclareMathOperator{\tr}{tr}
\DeclareMathOperator{\KL}{KL}

\DeclareMathOperator{\Law}{Law}
\renewcommand{\theenumi}{\roman{enumi}}

\endlocaldefs

\begin{document}

\begin{frontmatter}
\title{Persistence Pays: Square-root condition-number dependence for Metropolis-adjusted kinetic Langevin Monte Carlo with a non-reversibly updated accept/reject uniform}
\runtitle{Persistence Pays}

\begin{aug}
\author[A,B]{\fnms{Nawaf}~\snm{Bou-Rabee}}
\address[A]{Department of Mathematical Sciences, Rutgers University}
\address[B]{Center for Computational Mathematics, Flatiron Institute}
\end{aug}

\begin{abstract}
We prove accelerated mixing for a Metropolis-adjusted kinetic Langevin sampler: the OBABO splitting with friction $2\sqrt m$, adjusted with Neal's non-reversibly updated uniform, which is redrawn after every epoch of $E$ steps. For targets $e^{-U}$ on $\R^d$ with $m\Id\preceq\nabla^2U\preceq L\Id$ and an $L_H$-Lipschitz Hessian, started from a Gaussian at an approximate minimizer, the chain reaches total variation accuracy $\varepsilon$ after $\tilde O\big(\sqrt{\kappa d}+\sqrt{L_H}(d/m)^{3/4}\big)$ gradient evaluations, where $\kappa=L/m$ and $\tilde O$ hides logarithmic factors. The $\sqrt\kappa$ term is necessary, already on Gaussian targets. The proof discretizes Lu's entropic hypocoercivity by means of an entropy inequality for free transport, so that the relative entropy of the unadjusted chain decays at a rate of order $\sqrt m$ up to a discretization bias. With the non-reversibly updated uniform, the leading energy errors of velocity Verlet telescope along a path. A coupling argument combines the two. On Bayesian logistic regression posteriors for two standard data sets, which satisfy these assumptions, the gradient cost of the chain with the non-reversibly updated uniform grows roughly like the square root of the condition number at the posterior mode, and that of standard adjusted OBABO roughly like its $3/4$ power.
\end{abstract}


\end{frontmatter}

\section{Introduction}\label{sec:intro}

Let $U:\R^d\to\R$ satisfy $m\Id\preceq\nabla^2U\preceq L\Id$ and write $\kappa=L/m$. Kinetic Langevin dynamics with friction $\gamma=2\sqrt m$ converges to $\mu\propto e^{-U(x)-|v|^2/2}$ at a rate of order $\sqrt m$, independently of $L$: in $L^2(\mu)$ by Cao, Lu and Wang \cite{CLW23}, and in relative entropy by Lu \cite{Lu26}. Since stability of the explicit Verlet-based discretization requires a step of order $L^{-1/2}$, a discretization that inherited this rate would reach accuracy $\varepsilon$ after of order $\sqrt\kappa\log(1/\varepsilon)$ gradient evaluations, ignoring the dependence on $d$, against $\kappa\log(1/\varepsilon)$ for overdamped samplers. This is the sense of acceleration used here. Even for the class of centered Gaussian targets, sampling to sufficiently small fixed total variation error requires, in the worst case, at least a constant multiple of $\min\{\sqrt\kappa,d\}$ evaluations of $U$ and $\nabla U$ \cite[Corollary~4.4]{CDLLN23}. Thus, when $d\gtrsim\sqrt\kappa$, at least order $\sqrt\kappa$ such evaluations are necessary.

Whether the square-root acceleration of kinetic Langevin dynamics survives numerical discretization remains open for general strongly log-concave targets. The same question arises for numerical randomized HMC. We address it for Metropolis-adjusted OBABO under a Lipschitz-Hessian assumption and from a warm start.

\subsection{The sampler and the main result}

\begin{figure}[t]
\centering
\includegraphics[width=0.70\textwidth]{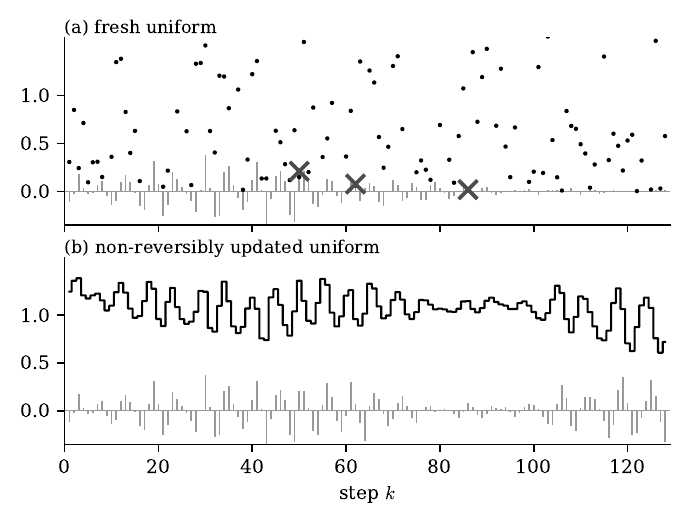}
\caption{Adjusted OBABO with (a) a fresh uniform and (b) a non-reversibly updated uniform; a move is accepted if $\Delta H<\ell$. Gray: energy error $\Delta H$. Black: threshold $\ell$, redrawn at every step in (a). Crosses: rejections.}
\label{fig:schematic}
\end{figure}

OBABO is the splitting of kinetic Langevin dynamics with step size $h>0$ that applies $O_{h/2}$, $B_{h/2}$, $A_h$, $B_{h/2}$, $O_{h/2}$ in this order, where $O$ is the exact Ornstein-Uhlenbeck flow of the velocity, $B$ a kick by the force $-\nabla U$ and $A$ a free flight \cite{BP07,LM13}. The composition $B_{h/2}\circ A_h\circ B_{h/2}$ of its three middle steps is the velocity Verlet integrator. We take $\gamma=2\sqrt m$ and write $\bar h=h\sqrt L$ for the dimensionless step size. The adjusted chain accepts or rejects the velocity Verlet proposal and reverses the velocity on rejection \cite{Hor91}. Instead of a fresh uniform at every step we use the non-reversibly updated uniform of Neal \cite{Nea20}. Written as $\ell=-\log u$, the standard test accepts a move with energy error $\Delta H$ when $\ell>\Delta H$, with a new $\ell\sim\Exp(1)$ at every step. In the variant analyzed here, the chain runs for $N$ epochs of $E$ consecutive steps each; $\ell$ is drawn once at the start of each epoch, and each accepted move applies Neal's update $\ell\leftarrow\ell-\Delta H$ (Definition~\ref{alg:pu}). Thus the value of the threshold before step $k+1$ of an epoch is its initial value minus the sum of the energy errors of the moves accepted among the first $k$ steps. Figure~\ref{fig:schematic} compares the two rules on one run, on the Sonar posterior of Section~\ref{sec:numerics}. The benefit is that a rejection within an epoch is then governed by the accumulated energy error, in which most of the per-step errors cancel, rather than by the sum of the per-step rejection probabilities, which grows with the number of steps (see Figure~\ref{fig:neal}). This chain leaves $\mu$ invariant (Lemma~\ref{lem:invariance}).

Informally, Theorem~\ref{thm:main} says the following. Assume that $\nabla^2U$ is $L_H$-Lipschitz and let $\nu$ be an initial law with $\|\dd\nu/\dd\mu\|_{L^2(\mu)}\le M_2$. Then total variation accuracy $\varepsilon$ is reached after
\begin{equation}\label{eq:intro-complexity}
O\Big(\Lambda^{3/2}\big(\sqrt{\kappa D}+\sqrt{l_H}\,D^{3/4}\big)\log(2/\varepsilon)\Big),\qquad D=d+\log(M_2/\varepsilon),
\end{equation}
gradient evaluations, where $l_H=L_H/m^{3/2}$ is the scale-invariant Hessian-Lipschitz ratio and $\Lambda=\log(2+\kappa+l_H+D)$. From a Gaussian started at an approximate minimizer, $\log M_2=O(d\log(2\kappa))$ (Remark~\ref{rem:warm}), and \eqref{eq:intro-complexity} becomes $\tilde O\big(\sqrt{\kappa d}+\sqrt{L_H}(d/m)^{3/4}\big)$, where $\tilde O$ hides logarithmic factors.

Up to logarithms, the first term of \eqref{eq:intro-complexity} dominates when $l_H\le\kappa/\sqrt d$ from this start; the exact threshold is $l_H\le\kappa/\sqrt D$. This condition does not require $U$ to be a small perturbation of a quadratic potential: it allows the Hessian to change by order $L$ over a distance $\sqrt{d/m}$, which is the largest scale of the position marginal, since $\E_{X\sim\mu_x}|X-x_*|^2\le d/m$, where $x_*$ is the minimizer of $U$ and $\mu_x$ the position marginal of $\mu$. The logistic regression posteriors of Section~\ref{sec:numerics} satisfy the assumptions of Theorem~\ref{thm:main}, but along that family $L_H$ grows like $\kappa^{3/2}$, so the second term dominates and the theorem gives cost of order $\kappa^{3/4}$ there; the measured cost grows like the square root of the condition number at the mode, beyond what the theorem proves.

\subsection{Ingredients of the proof}\label{sec:outline}

The proof has five ingredients, (a) to (e), described in turn.

\emph{(a) An entropy inequality for free transport.} Proposition~\ref{thm:entropy} states that for random vectors $X,V$ in $\R^d$ and any nondegenerate Gaussian density $g$, the quantity $\mathcal E_g(t)=\KL(\Law(X+tV,V)\,\|\,\Law(X+tV)\otimes g)$ satisfies $\mathcal E_g(t)\le\mathcal E_g(0)+\mathcal E_g(1)$ for $t\in[0,1]$. It follows from strong subadditivity of entropy applied to a Brownian bridge.

\emph{(b) Discrete entropic hypocoercivity.} Lu \cite{Lu26} proved accelerated entropy decay for kinetic Langevin dynamics with the functional $\Ent+\eta\,\Ccor$ for a suitable $\eta>0$, where $\Ent$ is relative entropy with respect to $\mu$ and $\Ccor(p)=\E_p\langle v,x-T(x)\rangle$, with $T$ the Brenier map from the position marginal of $p$ to $\mu_x\propto e^{-U}$. Its proof does not rely on an upper Hessian bound. Every step of it has an exact counterpart for OBABO except one: the force acts impulsively at the two ends of a free flight, and the resulting mismatch is the error of the trapezoidal rule applied to the differential entropy of the position marginal (Proposition~\ref{prop:DS0}). This error equals the excess of the velocity relative entropy along the free flight over its trapezoidal average, which (a) controls. The result is our second theorem, Theorem~\ref{thm:KL}: for unadjusted OBABO with $\sqrt mh\le3/16$ and every $U$ with $m\Id\preceq\nabla^2U\preceq L\Id$ and an $L_H$-Lipschitz Hessian, the relative entropy with respect to $\mu$ decays from any starting point at rate $\lambda\ge e^{-2\sqrt mh}\sqrt m/68$, which does not depend on $\kappa$, up to a bias term of order $\bar h^2$ times moments of the energy and of the velocity, plus a term proportional to $L_H$. With the moment bounds of (e), the bias is of order $\bar h^2D+l_H\bar h^2D^{3/2}/\kappa$, with $D$ as in \eqref{eq:intro-complexity}.

\emph{(c) The non-reversibly updated uniform.} Let $H(x,v)=U(x)-\min U+\frac12|v|^2$. For a velocity Verlet step from position $x_0$ to position $x_1$ with velocity $V$ during the free flight, the energy error is exactly $\Delta H=G(x_1)-G(x_0)+\tau$, with $G=\frac{h^2}8|\nabla U|^2$ and $\tau$ a trapezoidal rule error bounded by $\frac{L_H}{12}h^3|V|^3$ (Lemma~\ref{lem:energy}). The Ornstein-Uhlenbeck half-steps do not move positions, so along a path the $G$ terms telescope. With the non-reversibly updated uniform the chain accepts every move of an epoch exactly when the accumulated energy error stays below an exponential random variable (Proposition~\ref{prop:allaccept}). The rejection probability of an epoch is then of order $\bar h^2D$ plus the accumulated $\tau$, by the moment bounds of (e), with no factor of the number of steps.

\emph{(d) A coupling argument.} Run two copies of the adjusted chain, one from the initial distribution $\nu$ and one from the stationary law $\mu$, and drive them within each epoch by the Gaussian noises of the epoch, drawn from a maximal coupling of the two unadjusted chains (Lemma~\ref{lem:glue}). On the event that neither copy rejects, adjusted and unadjusted paths coincide, so the adjusted copies meet at the end of the epoch except with probability at most the total variation distance between the unadjusted endpoint laws, which (b) and Pinsker's inequality control, plus the two rejection probabilities. This is the localization idea of Bou-Rabee and Oberd\"orster \cite{BRO24}, in a form where the contraction is a total variation bound at the end of the epoch rather than a one-step Wasserstein contraction. The adjusted chain must start every epoch in a set $\mathcal D$ of points from which an epoch is accepted with high probability. Under $\mu$ each epoch starts with law $\mu$, so a union bound over the $N$ epochs of the run gives probability at most $N\mu(\mathcal D^c)$ that some epoch starts outside $\mathcal D$, and the Cauchy-Schwarz inequality transfers this bound to an $L^2$ warm start (Proposition~\ref{prop:coupling}). No bound on the excursions of the adjusted chain is needed.

\emph{(e) Moment bounds.} Theorem~\ref{thm:KL} needs moments of the unadjusted chain over one epoch, including third moments of the velocity. Below a velocity $v_*\asymp\kappa/(l_H\bar h^2)$ the energy error of velocity Verlet is dominated by friction, and a Lyapunov function, the energy with a small twist $\langle x,v\rangle$, gives exponential moments for the chain stopped when it leaves this region (Lemmas~\ref{lem:drift} and~\ref{lem:stopped}). Above $v_*$ we have no control of the unadjusted chain, uniformly over the class. We therefore apply the argument of (b) to a comparison chain that reverses the velocity instead of making a flight faster than $v_*$ (Section~\ref{sec:comparison}). It coincides with the unadjusted chain except on an event whose probability is at most $E$ times a quantity exponentially small in $v_*^2$. Each of its velocity Verlet steps multiplies $H-G$ by at most $1+0.081\sqrt m\,h$, and the noise adds at most $2(d+2)\sqrt m$ per unit time to the $L^2$ norm of $H-G$. The discrete hypocoercivity argument still applies to it, with an additional term controlled by the probability of a flight faster than $v_*$.

\subsection{Sketch of the proof of Theorem~\ref{thm:main}}

We refer to the ingredients (a) to (e) of Section~\ref{sec:outline}. The proof rests on the contraction guarantee of Theorem~\ref{thm:KL}: the relative entropy of the unadjusted chain decays at rate of order $\sqrt m$ in physical time, uniformly in $\kappa$, up to a discretization error. This is (b), and it is possible because of (a): every other step of Lu's argument transfers to OBABO exactly, and the one discrepancy, the trapezoidal error produced by the impulsive kicks, is bounded by (a). The remaining parts turn this contraction into a bound for the adjusted chain. Let $\mathcal D=\{H\le R\}$ with $R$ of order $D$, so that $N\mu(\mathcal D^c)$ is small (Lemma~\ref{lem:tail}). By (d), it suffices to show for every $z\in\mathcal D$ that (i) the unadjusted epoch from $z$ ends within a fixed small total variation distance of $\mu$, and (ii) the adjusted epoch from $z$ is rejection-free except with small probability; the two copies then meet at the end of each epoch with probability at least $3/4$, and accuracy $\varepsilon$ is reached after $N\asymp\log(1/\varepsilon)$ epochs. The moment bounds of (e) enter both. For (i), they bound the discretization error in Theorem~\ref{thm:KL}, applied to the comparison chain of (e), which coincides with the unadjusted chain except with negligible probability; the theorem then gives relative entropy at most a constant after a physical time of a fixed multiple of $\Lambda/\sqrt m$, that is, $E\asymp\Lambda\sqrt\kappa/\bar h$ steps, and Pinsker's inequality converts this to total variation. For (ii), the exact formula of (c) bounds the rejection probability by $\bar h^2$ times the maximal energy along the epoch plus the accumulated $\tau$, which these moments make of order $\bar h^2D$, small for $\bar h\asymp D^{-1/2}$. Thus $NE$ is \eqref{eq:intro-complexity}; the second term there comes from the smaller step forced by the $L_H$ terms in (b) and (e). Section~\ref{sec:proofmain} carries this proof out with explicit constants.

\subsection{The role of the Hessian-Lipschitz constant}\label{sec:LH}

The Hessian-Lipschitz assumption enters the proof at one point: the trapezoidal error $\tau$ of a velocity Verlet step is bounded by $\frac{L_H}{12}h^3|V|^3$ (Lemma~\ref{lem:energy}), and without such a bound the energy error of a step has no cancellation structure beyond the telescoping of $G$. This cubic bound is used twice. In (c) it keeps the accumulated energy error of an epoch, and hence the rejection probability, small; in (b) it is the term proportional to $L_H$ in Theorem~\ref{thm:KL}. Both uses require the flight speed to stay below $v_*\asymp\kappa/(l_H\bar h^2)$, where friction dominates the error, which is what the comparison chain of (e) enforces. The constant $L_H$ then appears in Theorem~\ref{thm:main} only through the second term of \eqref{eq:intro-complexity} and the corresponding restriction on the step. This assumption cannot be dropped from the present argument: with only strong convexity and a Lipschitz gradient, the available bound on $\tau$, summed over an epoch, is not small at the step size of Theorem~\ref{thm:main}, and exploiting cancellation of the signed errors would require a different argument (Section~\ref{sec:discussion}).

\subsection{Relation to the literature}

\emph{Metropolis-adjusted kinetic samplers.} For an adjusted OABAO splitting, Bou-Rabee and Oberd\"orster \cite{BRO24} proved mixing after $\tilde O\big(\kappa^{3/2}\max\{L_H^{1/2}(d/m)^{3/4},(\kappa d)^{1/2}\}\log^2(1/\varepsilon)\big)$ gradient evaluations from cold starts with initial energy $H(z)=O(\kappa d)$, and noted that the optimal dependence on $\kappa$ was unknown. Apers, Gribling and Szil\'agyi \cite{AGS24} obtained $\tilde O(\sqrt\kappa\,d^{1/4}\log(1/\varepsilon))$ gradient evaluations on Gaussian targets, using an unadjusted warm-up followed by adjusted HMC with random integration times; their proof relies on Gaussian structure. The general-class upper bounds in \cite{LST20,WSC22} for MALA, equivalently HMC with one leapfrog step, depend linearly on $\kappa$; related lower bounds for MALA and HMC with fixed integration parameters appear in \cite{LST21}. To our knowledge, Theorem~\ref{thm:main} is the first bound of order $\sqrt\kappa$ on the number of gradient evaluations for a Metropolis-adjusted discretization of kinetic Langevin dynamics on a class of non-Gaussian targets. Neal \cite{Nea20} introduced the non-reversibly updated uniform and established its validity; its computational benefits have been investigated empirically \cite{Nea20,HS22}. To our knowledge, this is the first nonasymptotic mixing bound for a kinetic Langevin discretization using this update.

\emph{Accelerated samplers on the full class.} From a $\chi^2$ warm start, the proximal bouncy particle sampler of Chen, Chewi, Lu and Zhang \cite{CCLZ26} and the reflected proximal dynamics of Yuan, Guo and Chen \cite{YGC26} achieve total variation error $\varepsilon$ using $\tilde O(\sqrt\kappa\,d^{1/4})$ expected oracle queries, with polylogarithmic dependence on $1/\varepsilon$. Both use event-driven proximal constructions, controlled approximations of conditional sampling and event simulation, and $L^2$ contraction estimates. Mitra, Srinivasan, Wang and Wibisono \cite{MSWW26} proved accelerated mixing for randomized Hamiltonian Monte Carlo with the exact flow. Lu and Luo \cite{LL26} combined continuous-time entropic hypocoercivity with exact simulation by windowed thinning to obtain $\tilde O(d\sqrt\kappa)$ expected gradient evaluations from a Gaussian start centered at a supplied minimizer. The proximal methods require no Hessian-Lipschitz assumption, and their $d^{1/4}$ dependence improves on the $d^{1/2}$ factor in our $\sqrt\kappa$ term. Our contribution is an acceleration guarantee for a Metropolis-adjusted discretization.

\emph{Unadjusted discretizations.} Altschuler, Chewi and Zhang \cite{ACZ25} analyzed a variant of the randomized midpoint method of Shen and Lee \cite{SL19}, obtaining $\tilde O(\kappa^{5/6}d^{5/3}\varepsilon^{-2/3})$ gradient evaluations to reach total variation error $\varepsilon$ from a warm start, and exact simulation of the diffusion \cite{CCRZ26} reaches $\kappa^{2/3}d^{1/3}$ up to logarithms in $1/\varepsilon$. Bou-Rabee \cite{BR26a} proved that standard Strang splittings, including OBABO and BAOAB, are not accelerated from cold starts in total variation when the step size and friction are chosen from $(m,L,d)$ alone. Eliminating the velocity yields a noisy heavy-ball recursion: either a Gaussian target relaxes slowly, or attracting cycles arising from the construction of Goujaud, Taylor and Dieuleveut \cite{GTD25} produce metastability under noise. Our result concerns warm starts and a Metropolis-adjusted chain, to which these lower bounds do not apply.

\emph{Discrete hypocoercivity.} Camrud, Durmus, Monmarch\'e and Stoltz \cite{CDMS23} proved discrete-time entropic hypocoercivity for unadjusted generalized HMC using a modified entropy of Villani type \cite{Vil09}. Monmarch\'e and Schuh \cite{MS24} extended this approach to an unadjusted second-order splitting scheme for mean-field kinetic Langevin dynamics. Theorem~\ref{thm:KL} gives an analogous bound for unadjusted OBABO, using entropy augmented by Lu's transport corrector \cite{Lu26}, with a decay rate of order $\sqrt m$ in physical time, uniformly in the condition number, up to a discretization error that depends on moments of the chain.

\subsection{Organization}

Section~\ref{sec:setting} sets notation, defines the chain and states the main results. Section~\ref{sec:entropy} proves the entropy inequality for free transport. Section~\ref{sec:hypo} develops the discrete entropic hypocoercivity and proves Theorem~\ref{thm:KL}. Section~\ref{sec:coupling} characterizes the event that every move of an epoch is accepted, and contains the coupling argument, the moment bounds and the comparison chain. Section~\ref{sec:proofmain} proves Theorem~\ref{thm:main}, through explicit conditions on the parameters (Proposition~\ref{prop:explicit}). Section~\ref{sec:numerics} applies the results to Bayesian logistic regression and compares the two acceptance rules, and Section~\ref{sec:discussion} discusses open problems. Appendix~\ref{app:free} proves the free-transport estimate used in Lemma~\ref{lem:approx}, and Appendix~\ref{app:sharp} proves Proposition~\ref{prop:sharp}.

\subsection{Formal verification}

The statements of Theorems~\ref{thm:main} and~\ref{thm:KL} and Propositions~\ref{prop:sharp} and~\ref{thm:entropy}, with the definitions they rest on, have been formalized and their proofs checked in Lean~4 against the Mathlib library, with no axioms beyond the standard three of Lean's kernel. The formalization follows the proofs given here, including the explicit constants, and in places proves slightly more (Proposition~\ref{prop:sharp} for arbitrary friction and at every step within an epoch; Proposition~\ref{thm:entropy} without moment assumptions). The Lean sources, a table matching each formal statement to its counterpart in this paper, and the build log are available at \url{https://github.com/nawafbourabee/persistence-pays-lean} (archived at doi:10.5281/zenodo.23671351).

\section{Setting and main results}\label{sec:setting}

\subsection{Notation and assumptions}

Points of phase space are $z=(x,v)\in\R^d\times\R^d$. We write $\langle\cdot,\cdot\rangle$ and $|\cdot|$ for the Euclidean inner product and norm and, for a matrix $A\in\R^{d\times d}$, $\|A\|=\sup_{|u|=1}|Au|$ for its operator norm. We assume:

\begin{itemize}
\item[(A1)] $U\in C^2(\R^d)$ and $m\Id\preceq\nabla^2U(x)\preceq L\Id$ for all $x$, with $0<m\le L$.
\item[(A2)] $\|\nabla^2U(x)-\nabla^2U(y)\|\le L_H|x-y|$ for all $x,y$.
\end{itemize}

Let $x_*$ be the minimizer of $U$ and $U_*=U(x_*)$. Put $H(x,v)=U(x)-U_*+\frac12|v|^2$, $\mu(\dd z)\propto e^{-H(z)}\dd z$, $\mu_x(\dd x)\propto e^{-U(x)}\dd x$, and let $\varphi$ be the standard Gaussian law on $\R^d$, so $\mu=\mu_x\otimes\varphi$. We use the dimensionless quantities
\[
\kappa=L/m,\qquad \bar h=h\sqrt L,\qquad l_H=L_Hm^{-3/2},\qquad\text{so that}\qquad\sqrt m\,h=\bar h/\sqrt\kappa .
\]
For probability measures $p,q$, $\KL(p\|q)$ is the relative entropy and $\mathrm{TV}(p,q)=\sup_A|p(A)-q(A)|$ is the total variation distance. For a law $p$ and a measurable map or Markov kernel $Q$, $pQ$ denotes the law after applying $Q$. $\Exp(1)$ is the standard exponential law. The differential entropy of a law $\rho$ with density on $\R^n$ is $\mathrm h(\rho)=-\int\rho\log\rho$; for a random vector $Y$ we write $\mathrm h(Y)=\mathrm h(\Law(Y))$.

\subsection{OBABO}

Let $h>0$ be the step size. For $s\ge0$ and $\xi\in\R^d$ define the maps
\begin{gather*}
O_s^\xi(x,v)=\big(x,\,e^{-\gamma s}v+(1-e^{-2\gamma s})^{1/2}\xi\big),\\
B_s(x,v)=(x,v-s\nabla U(x)),\qquad A_s(x,v)=(x+sv,v),
\end{gather*}
and let $O_s$ be the Markov kernel $O_s(z,\cdot)=\Law(O_s^\xi(z))$ with $\xi\sim\varphi$. We call $O_{h/2}^{\xi}$ an Ornstein-Uhlenbeck half-step. Throughout, $\gamma=2\sqrt m$, $r=e^{-\gamma h/2}=e^{-\sqrt mh}$ and $\sigma=(1-r^2)^{1/2}$. Velocity Verlet is $\Phi=B_{h/2}\circ A_h\circ B_{h/2}$, and the OBABO kernel is $K(z,\cdot)=\Law\big((O_{h/2}^{\xi'}\circ\Phi\circ O_{h/2}^{\xi})(z)\big)$ where $\xi$ and $\xi'$ are independent random vectors, each with law $\varphi$. For $w=(x,v)$ write
\begin{equation}\label{eq:halfvel}
V(w)=v-\tfrac h2\nabla U(x),\qquad \Phi(w)=\big(x_1,\,V(w)-\tfrac h2\nabla U(x_1)\big),\quad x_1=x+hV(w),
\end{equation}
and define
\[
G(x)=\tfrac{h^2}{8}|\nabla U(x)|^2,\qquad \tau(w)=U(x_1)-U(x)-\tfrac12\langle\nabla U(x)+\nabla U(x_1),x_1-x\rangle .
\]
Let $\xi_0,\xi_0',\xi_1,\xi_1',\dots$ be independent random vectors in $\R^d$, each with law $\varphi$. The unadjusted chain from $z$ is $Z_0=z$, $W_k=O_{h/2}^{\xi_k}(Z_k)$, $Z_{k+1}=O_{h/2}^{\xi_k'}(\Phi(W_k))$. Let $X_k$ be the position of $Z_k$ and $V_k=V(W_k)$. Thus $W_k$ also has position $X_k$, $V_k$ is the velocity during the free flight from $X_k$ to $X_{k+1}=X_k+hV_k$, and $\Phi(W_k)$ is the result of the $(k+1)$-th velocity Verlet step.

\begin{lemma}[Energy error of velocity Verlet]\label{lem:energy}
Under (A1), for every $w$,
\begin{gather*}
\Delta H(w):=H(\Phi(w))-H(w)=G(x_1)-G(x)+\tau(w),\\
|\tau(w)|\le\tfrac{\bar h^2}{8}|V(w)|^2,\qquad 0\le G\le\tfrac{\bar h^2}4(U-U_*).
\end{gather*}
Under (A2) moreover $|\tau(w)|\le\frac{L_H}{12}h^3|V(w)|^3$, and $\tau\equiv0$ when $U$ is quadratic.
\end{lemma}

\begin{proof}
Write $g_0=\nabla U(x)$, $g_1=\nabla U(x_1)$ and $V=V(w)$. Then $v=V+\frac h2g_0$ and the final velocity is $V-\frac h2g_1$, so
\[
\tfrac12|V-\tfrac h2g_1|^2-\tfrac12|V+\tfrac h2g_0|^2=-\tfrac h2\langle V,g_0+g_1\rangle+\tfrac{h^2}8(|g_1|^2-|g_0|^2).
\]
Since $hV=x_1-x$, adding $U(x_1)-U(x)$ gives the identity. Put $f(t)=U(x+t(x_1-x))$. Integration by parts gives $\tau=-\int_0^1(t-\frac12)f''(t)\dd t$. Thus $\tau=\int_0^{1/2}(\frac12-t)f''(t)\dd t-\int_{1/2}^1(t-\frac12)f''(t)\dd t$ is the difference of two numbers in $[0,\frac L8|x_1-x|^2]$, since by (A1) $0\le f''\le L|x_1-x|^2$, so $|\tau|\le\frac L8|x_1-x|^2=\frac{\bar h^2}8|V|^2$. Since $\int_0^1(t-\frac12)\dd t=0$, also $\tau=-\int_0^1(t-\frac12)(f''(t)-f''(\frac12))\dd t$, and under (A2) $|f''(t)-f''(\frac12)|\le L_H|t-\frac12||x_1-x|^3$, which gives $|\tau|\le\frac{L_H}{12}|x_1-x|^3$. Finally, the descent lemma gives $U_*\le U(x-\frac1L\nabla U(x))\le U(x)-\frac1{2L}|\nabla U(x)|^2$, so $|\nabla U|^2\le2L(U-U_*)$ and $0\le G\le\frac{\bar h^2}4(U-U_*)$.
\end{proof}

\subsection{The adjusted chain}

\begin{definition}[Adjusted OBABO with a non-reversibly updated uniform]\label{alg:pu}
Fix $E\in\mathbb N$ and let $S(x,v)=(x,-v)$ be the velocity reversal. One epoch from $z\in\R^{2d}$ is:
\begin{enumerate}
\item draw $\ell\sim\Exp(1)$, independent of the past;
\item for $k=0,\dots,E-1$, with $w=O_{h/2}^{\xi_k}(z)$: set $(z,\ell)\leftarrow\big(O_{h/2}^{\xi_k'}(\Phi(w)),\,\ell-\Delta H(w)\big)$ if $\ell>\Delta H(w)$ (accept), and $(z,\ell)\leftarrow\big(O_{h/2}^{\xi_k'}(S(w)),\,\ell\big)$ otherwise (reject).
\end{enumerate}
Here $\ell,\xi_0,\xi_0',\dots,\xi_{E-1},\xi_{E-1}'$ are independent, and each $\xi_k$ and $\xi_k'$ has law $\varphi$, as in the unadjusted chain. We write $P$ for the Markov kernel of one epoch on $\R^{2d}$.
\end{definition}

At the start of each epoch, $u=e^{-\ell}$ is uniform on $(0,1)$. In terms of $u$, the test $\ell>\Delta H$ is the Metropolis test $u<e^{-\Delta H}$, and the update $\ell\leftarrow\ell-\Delta H$ is Neal's update $u\leftarrow ue^{\Delta H}$ \cite[Sections~2 and~4]{Nea20}. Neal also moves $u$ by a small amount between iterations; here a single uniform is drawn at the start of each epoch, and its running value changes only through the acceptance update $u\leftarrow ue^{\Delta H}$, remaining unchanged on rejection. After one initial gradient evaluation, each step uses one new gradient evaluation. On acceptance, the gradient at the proposed position is retained; on rejection, the gradient at the current position is retained.

\begin{figure}[t]
\centering
\includegraphics[width=0.7\textwidth]{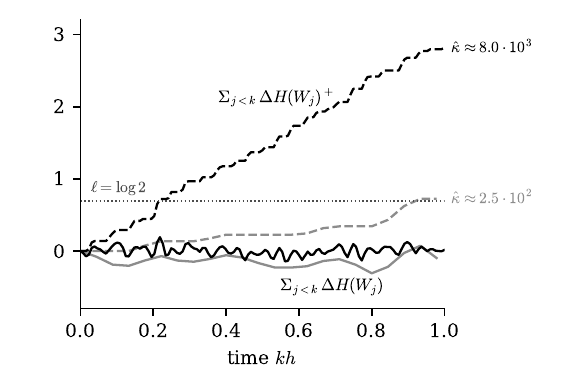}
\caption{Accumulated energy error (solid) and accumulated positive parts (dashed) along unadjusted OBABO on the Sonar posterior of Section~\ref{sec:numerics}. Dotted: the median of $\ell$.}
\label{fig:neal}
\end{figure}

As long as every move of an epoch is accepted, the threshold before step $k+1$ is $\ell-\sum_{j<k}\Delta H(W_j)$, where $\ell$ is the value drawn in (i). So every move of the epoch is accepted exactly when the accumulated energy error $\sum_{j<k}\Delta H(W_j)$ stays below $\ell$ for all $k\le E$ (Proposition~\ref{prop:allaccept}). By Lemma~\ref{lem:energy} the terms $G(X_{j+1})-G(X_j)$ telescope in this sum. With a fresh uniform at every step, the probability that some move of the epoch is rejected is bounded by $\sum_{j<E}\E[1-e^{-\Delta H(W_j)^+}]$, and the positive parts do not cancel. Figure~\ref{fig:neal} shows both sums along one unit of time.

\begin{lemma}[Invariance]\label{lem:invariance}
Each step of (ii) in Definition~\ref{alg:pu} preserves $\mu\otimes\Exp(1)$ on $\R^{2d}\times(0,\infty)$, and $\mu P=\mu$.
\end{lemma}

\begin{proof}
We apply the auxiliary-variable involution framework of \cite{BCM26}, with the threshold $\ell$ as the auxiliary variable, on the extended state space $\mathbb S=\R^{2d}\times(0,\infty)$ with target $\bar\mu=\mu\otimes\Exp(1)$. The map $\Phi$ is a composition of symplectic maps, so it preserves Lebesgue measure, and $\Phi^{-1}=S\circ\Phi\circ S$. Hence $S\circ\Phi$ is a volume-preserving involution, and $H\circ S=H$. Let $\mathcal A=\{(w,\ell)\in\mathbb S:\ell>\Delta H(w)\}$, and define $\Upsilon$ on $\mathbb S$ by
\[
\Upsilon(w,\ell)=((S\circ\Phi)(w),\ell-\Delta H(w))\ \text{on }\mathcal A,\qquad\Upsilon=\mathrm{id}\ \text{on }\mathbb S\setminus\mathcal A .
\]
The two branches of $\Upsilon$ correspond to the accept/reject step in (ii): $\mathcal A$ is the event that the Verlet proposal is accepted, and off $\mathcal A$ the map $\Upsilon$ leaves the state unchanged; composing with the extended velocity reversal $\hat S$ defined at the end of the proof gives the rejection update $w\mapsto S(w)$ with $\ell$ fixed.

We verify that $\Upsilon$ is a volume-preserving involution of $\mathbb S$, as required by the involution framework of \cite{BCM26}, and that it preserves the extended energy.
\begin{itemize}
\item \emph{$\Upsilon$ preserves the acceptance region.} Let $(w,\ell)\in\mathcal A$ and $(w',\ell')=\Upsilon(w,\ell)$. Then $\ell'=\ell-\Delta H(w)>0$. Since $H\circ S=H$ and $(S\circ\Phi)\circ(S\circ\Phi)=\mathrm{id}$, we have $\Delta H(w')=H((\Phi\circ S\circ\Phi)(w))-H((S\circ\Phi)(w))=H(w)-H(w')=-\Delta H(w)$, so $\ell'>\Delta H(w')$ is equivalent to $\ell>0$, which holds on $\mathbb S$. Hence $\Upsilon(\mathcal A)\subseteq\mathcal A$, and $\Upsilon$ maps $\mathbb S$ into itself.
\item \emph{Involution.} For $(w,\ell)\in\mathcal A$, $\Upsilon(w',\ell')=((S\circ\Phi\circ S\circ\Phi)(w),\ell'+\Delta H(w))=(w,\ell)$, so $\Upsilon\circ\Upsilon=\mathrm{id}$ on $\mathcal A$. On the complement this is trivial. Thus $\Upsilon$ is an involution of $\mathbb S$.
\item \emph{Jacobian.} The set $\mathcal A=\{(w,\ell)\in\mathbb S:\ell-\Delta H(w)>0\}$ is open, since $\Delta H$ is continuous. Since $U\in C^2$, on $\mathcal A$ the map $\Upsilon$ is $C^1$, and its Jacobian is block lower triangular with diagonal blocks $D(S\circ\Phi)$ and $1$. Hence $|\det D\Upsilon|=1$ on $\mathcal A$. Since $\Upsilon$ maps $\mathcal A$ bijectively onto itself, change of variables shows that it preserves Lebesgue measure on $\mathcal A$. It is the identity on the complement, so it preserves Lebesgue measure on $\mathbb S$.
\item \emph{Extended energy.} On $\mathcal A$, $H(w')+\ell'=H(w)+\Delta H(w)+\ell-\Delta H(w)=H(w)+\ell$; off $\mathcal A$ this holds trivially.
\end{itemize}
In this representation, the Metropolis kernel built from the volume-preserving involution $\Upsilon$ proposes $\Upsilon(w,\ell)$ and accepts it with probability
\[
\min\Big\{1,\ \frac{e^{-H(w')-\ell'}}{e^{-H(w)-\ell}}\Big\}=1,
\]
by the extended-energy identity. So this kernel reduces to the deterministic map $\Upsilon$, and by \cite[Lemma~2.2]{BCM26}, $\Upsilon$ preserves $\bar\mu$. The unit acceptance probability refers only to this extended-space representation. It does not mean that every Verlet proposal is accepted: rejection of a proposal corresponds to the branch $\mathbb S\setminus\mathcal A$, on which $\Upsilon$ is the identity.

Let $\hat S(w,\ell)=(S(w),\ell)$ on $\mathbb S$. The composition $\hat S\circ\Upsilon$ equals $\Phi(w)$ with $\ell\leftarrow\ell-\Delta H(w)$ on $\mathcal A$ and $S(w)$ with $\ell$ fixed off $\mathcal A$, which is exactly the accept/reject update of (ii). Thus one step of (ii) is $(w,\ell)\mapsto(O_{h/2}^{\xi_k}(w),\ell)$, followed by $\hat S\circ\Upsilon$, followed by $(w,\ell)\mapsto(O_{h/2}^{\xi_k'}(w),\ell)$. The maps $\Upsilon$ and $\hat S$ preserve $\bar\mu$, by the preceding argument and the identity $H\circ S=H$, respectively. The two Ornstein-Uhlenbeck half-steps, viewed as Markov kernels, also preserve $\bar\mu$: they preserve $\varphi$ and leave $x$ and $\ell$ unchanged. Hence each step of (ii) preserves $\bar\mu$. If an epoch starts from $z\sim\mu$, then after (i) the pair $(z,\ell)$ has law $\bar\mu$, since $\ell$ is drawn independently of $z$, and every step of (ii) preserves $\bar\mu$. So the final position and velocity have law $\mu$, that is, $\mu P=\mu$. This is the extended-space invariance argument of \cite{BCM26}; since $\ell$ is carried from step to step within an epoch rather than refreshed, we do not use the reversibility statement of \cite[Theorem~2.1]{BCM26} for a fully refreshed auxiliary variable.
\end{proof}

\subsection{Main results}

Theorem~\ref{thm:main} bounds the number of gradient evaluations needed by the adjusted chain, and Proposition~\ref{prop:sharp} shows that its $\sqrt\kappa$ dependence is necessary. Theorem~\ref{thm:KL}, stated after them, is the decay of relative entropy for unadjusted OBABO on which the proof of Theorem~\ref{thm:main} rests. The rate below comes from Proposition~\ref{prop:R} with $\theta=r/16$:
\begin{equation}\label{eq:lambda}
\lambda=\frac{r^2\sqrt m}{4(16+r)}\ \ge\ \frac{\sqrt m}{68}\,e^{-2\sqrt mh},
\end{equation}
where $r=e^{-\sqrt mh}$ is the velocity contraction factor of an Ornstein-Uhlenbeck half-step. Let $t_E=\sqrt m\,Eh$ be the duration of an epoch in units of $1/\sqrt m$.

\begin{mainthm}[Mixing of the adjusted chain]\label{thm:main}
Assume (A1) and (A2). Let $\varepsilon\in(0,1)$ and let $\nu$ be a probability measure on $\R^{2d}$ with $M_2:=\|\dd\nu/\dd\mu\|_{L^2(\mu)}<\infty$. Put
\[
D=d+\log(M_2/\varepsilon),\qquad\Lambda=\log(2+\kappa+l_H+D),\qquad N=\lceil\log(2/\varepsilon)/\log4\rceil.
\]
There is an absolute constant $C_0$, and there are $\bar h$ and $E$ depending only on $\kappa$, $l_H$, $d$, $M_2$ and $\varepsilon$, such that
\[
NE\le C_0\,\Lambda^{3/2}\big(\sqrt{\kappa D}+\sqrt{l_H}\,D^{3/4}\big)\log(2/\varepsilon)
\]
and $\mathrm{TV}(\nu P^N,\mu)\le\varepsilon$. Thus accuracy $\varepsilon$ is reached after at most $NE+1$ gradient evaluations.
\end{mainthm}

The proof (Section~\ref{sec:proofmain}) takes $t_0=C\Lambda$ for a large absolute constant $C$, and
\[
\bar h^2=\min\Big\{\frac{1}{C^2Q\Lambda},\ \frac{\kappa}{C^2l_HQ^{3/2}t_0}\Big\},\qquad E=\Big\lceil\frac{t_0\sqrt\kappa}{\bar h}\Big\rceil+1,
\]
where $Q$, defined in Proposition~\ref{prop:explicit}, satisfies $D\le Q\le94D$. When $L_H=0$, the second entry in the minimum is interpreted as $+\infty$. Thus, up to logarithms, $\bar h$ is the smaller of $D^{-1/2}$ and $(\kappa/(l_HD^{3/2}))^{1/2}$, and an epoch lasts $t_E\asymp\Lambda$ units of $1/\sqrt m$. Proposition~\ref{prop:explicit} gives explicit conditions on $\bar h$ and $E$ under which $\mathrm{TV}(\nu P^N,\mu)\le\varepsilon$.

The $\sqrt\kappa$ dependence in Theorem~\ref{thm:main} is necessary for this algorithm, already on Gaussian targets and from a warm start.

\begin{proposition}[Sharpness in $\kappa$]\label{prop:sharp}
Let $d\ge2$, $0<m\le L$, $\kappa=L/m$ and $U(x)=\frac m2x_1^2+\frac L2\sum_{i=2}^dx_i^2$. Consider adjusted OBABO with friction $\gamma>0$ and step size $h$ such that $\bar h=h\sqrt L<2$, using either a fresh uniform at every step or the non-reversibly updated uniform of Definition~\ref{alg:pu} with any $E\in\mathbb N$. Start from $\mu_A=\mu(\cdot\,|\,x_1\le0)$, so that $\dd\mu_A/\dd\mu\le2$, and let $X_n$ be the position after $n$ steps, rejected ones included. Then $\mathrm{TV}(\Law(X_n),\mu_x)>\frac14$ whenever
\[
n<\frac{\sqrt{2\kappa}}{4\bar h}\Big(1+\frac{\bar h^6}{64}\Big)^{(d-1)/2},
\]
and hence, for every $\bar h\in(0,2)$, whenever $n<\frac{\sqrt2}{8}\sqrt\kappa\,d^{1/6}$.
\end{proposition}

The second condition follows from the first by minimizing over $\bar h$. Since $n$ steps cost $n+1$ gradient evaluations, for every fixed $d\ge2$ the factor $\sqrt\kappa$ in Theorem~\ref{thm:main} cannot be improved beyond logarithmic factors, and a warm start does not remove it. In the dimension, at the step size of Theorem~\ref{thm:main}, which on this Gaussian is $\bar h\asymp(D\Lambda)^{-1/2}$ since $L_H=0$, the first condition requires $n$ of order at least $\sqrt{\kappa D\Lambda}$, so there the bound $\sqrt{\kappa D}$ of Theorem~\ref{thm:main} is sharp up to logarithmic factors. On this Gaussian family, over all step sizes, the dimension dependence of this algorithm lies between $d^{1/6}$ and $d^{1/2}$; improving on $d^{1/2}$ would require step sizes larger than those of Theorem~\ref{thm:main}. The proof of Proposition~\ref{prop:sharp}, in Appendix~\ref{app:sharp}, combines a conductance bound for the half-space $\{x_1\le0\}$ with an exact Gaussian integral bounding the accepted crossing probability.

\begin{remark}[Gaussian starts]\label{rem:warm}
A Gaussian start centered at an approximate minimizer of $U$ is warm. Let $\hat x$ be a point with $L|\hat x-x_*|^2\le1$, that is, within distance $L^{-1/2}$ of the minimizer $x_*$, and let $\nu=\mathcal N(\hat x,(2L)^{-1}\Id)\otimes\varphi$. Then $M_2\le e(2\kappa)^{d/2}$. Indeed, since $e^{-(U-U_*)}\ge e^{-L|x-x_*|^2/2}$, $\int e^{-(U-U_*)}\le(2\pi/m)^{d/2}$ and $\frac12|x-x_*|^2\le|x-\hat x|^2+|\hat x-x_*|^2$,
\[
\frac{\dd\nu}{\dd\mu}(x,v)\le\Big(\frac{L}{\pi}\Big)^{d/2}\Big(\frac{2\pi}{m}\Big)^{d/2}e^{L|x-x_*|^2/2-L|x-\hat x|^2}\le(2\kappa)^{d/2}e^{L|\hat x-x_*|^2}.
\]
From a point $x_0$, Nesterov's accelerated gradient method \cite[Section~2.2]{Nes18} finds such an $\hat x$ after $O\big(\sqrt\kappa\log(2+\kappa L|x_0-x_*|^2)\big)$ gradient evaluations. Then $D\le1+d+\frac d2\log(2\kappa)+\log(1/\varepsilon)$, and Theorem~\ref{thm:main} gives $\tilde O\big(\sqrt{\kappa d}+\sqrt{L_H}(d/m)^{3/4}\big)$ gradient evaluations in addition to those of this method, using $\sqrt{l_H}\,d^{3/4}=\sqrt{L_H}(d/m)^{3/4}$.
\end{remark}

The second result is decay of relative entropy for the unadjusted chain.

\begin{mainthm}[Relative entropy decay of unadjusted OBABO]\label{thm:KL}
Assume (A1) and (A2), $\gamma=2\sqrt m$ and $\sqrt mh\le3/16$, and let $\lambda$ be given by \eqref{eq:lambda}. Then for every $z\in\R^{2d}$ and $n\ge1$,
\begin{multline*}
\KL(\delta_zK^{n+1}\|\mu)\le\tfrac{17}{15}e^{-\lambda nh}\Big[\E_zH(Z_1)+d\log\frac{\kappa}{2e\bar h^2(1-\sqrt mh)}\Big]\\
+\tfrac{16}{15}\Big[\tfrac{\bar h^2}{4}\E_z\big(U(X_{n+1})-U_*\big)+\sup_{1\le k\le n}\Big(\tfrac{L_Hh^2}{12\lambda}\E_z|V_k|^3+\tfrac{221}{16r}\bar h^2\E_z|V_k|^2\Big)\Big].
\end{multline*}
\end{mainthm}

The rate involves no $\kappa$. The upper Hessian bound enters only through the factor $\bar h^2$ of the bias term, and $L_H$ only through the term in $\E|V_k|^3$. Theorem~\ref{thm:KL} holds for every $U$ in the class, but its right side involves moments of the unadjusted chain, which we do not control uniformly over the class at flight speeds where friction no longer dominates the energy error of velocity Verlet (above the threshold $v_*$ of Section~\ref{sec:moments}). The proof of Theorem~\ref{thm:main} therefore applies the same argument to a comparison chain (Proposition~\ref{prop:KLsharp}), whose moments are controlled.

\section{An entropy inequality for free transport}\label{sec:entropy}

\begin{proposition}\label{thm:entropy}
Let $X,V$ be random vectors in $\R^d$, let $g$ be a nondegenerate Gaussian density on $\R^d$, and put $\mathcal E_g(t)=\KL\big(\Law(X+tV,V)\,\big\|\,\Law(X+tV)\otimes g\big)$. If $\mathcal E_g(0)$ and $\mathcal E_g(1)$ are finite, then
\[
\mathcal E_g(t)\le\mathcal E_g(0)+\mathcal E_g(1)\qquad\text{for all }t\in[0,1].
\]
\end{proposition}

We give two proofs. The first, which shows the mechanism, uses strong subadditivity of entropy along a Brownian bridge and assumes finite second moments and finite differential entropy. The second is a dual argument requiring only the finiteness of $\mathcal E_g(0)$ and $\mathcal E_g(1)$.

\subsection{Proof by strong subadditivity}

Write $Y=X+V$, $M_t=(1-t)X+tY=X+tV$, and let $g$ be the density of $\mathcal N(b,\Sigma)$. If $(X,V)$ has a density, the shears $(x,v)\mapsto(x+tv,v)$ have unit Jacobian, so $\mathrm h(M_t,V)=\mathrm h(X,Y)$ and
\begin{equation}\label{eq:Eg-entropy}
\mathcal E_g(t)=-\mathrm h(X,Y)+\mathrm h(M_t)-\E\log g(Y-X)
\end{equation}
whenever the three terms are finite. Hence $\mathcal E_g(t)\le\mathcal E_g(0)+\mathcal E_g(1)$ is equivalent to
\begin{equation}\label{eq:linent}
\mathrm h(X,Y)+\mathrm h(M_t)\le\mathrm h(X)+\mathrm h(Y)-\E\log g(Y-X).
\end{equation}

\begin{proof}[First proof of Proposition~\ref{thm:entropy}]
Assume that $(X,V)$ has finite second moments and finite entropy. Then all entropies below are finite: they are bounded above by those of Gaussians with the same second moments, and below by $\mathrm h(A)\ge\mathrm h(A,B)-\mathrm h(B)$. Fix $t\in(0,1)$, let $\zeta\sim\mathcal N(0,\Sigma)$ be independent of $(X,Y)$, and put $W=M_t+\sqrt{t(1-t)}\,\zeta$. Strong subadditivity of differential entropy gives
\[
\mathrm h(X,Y,W)+\mathrm h(W)\le\mathrm h(X,W)+\mathrm h(Y,W).
\]
Here $\mathrm h(X,Y,W)=\mathrm h(X,Y)+\frac d2\log(2\pi et(1-t))+\frac12\log\det\Sigma$ and $\mathrm h(W)\ge\mathrm h(M_t)$. Given $X$, $W-X=t(Y-X)+\sqrt{t(1-t)}\zeta$, so by Gibbs' inequality against the density of $\mathcal N(tb,t\Sigma)$,
\[
\mathrm h(W|X)\le\tfrac d2\log(2\pi t)+\tfrac12\log\det\Sigma+\tfrac t2\E|Y-X-b|_\Sigma^2+\tfrac{(1-t)d}2,
\]
with $|u|_\Sigma^2=u^\top\Sigma^{-1}u$. In the same way, against $\mathcal N(-(1-t)b,(1-t)\Sigma)$,
\[
\mathrm h(W|Y)\le\tfrac d2\log(2\pi(1-t))+\tfrac12\log\det\Sigma+\tfrac{1-t}2\E|Y-X-b|_\Sigma^2+\tfrac{td}2.
\]
Insert $\mathrm h(X,W)=\mathrm h(X)+\mathrm h(W|X)$ and $\mathrm h(Y,W)=\mathrm h(Y)+\mathrm h(W|Y)$. The terms in $\log t(1-t)$ cancel, and what remains is $\mathrm h(X,Y)+\mathrm h(M_t)\le\mathrm h(X)+\mathrm h(Y)+\frac d2\log(2\pi)+\frac12\log\det\Sigma+\frac12\E|Y-X-b|^2_\Sigma$, which is \eqref{eq:linent}.
\end{proof}

\subsection{The general case}

Without moment assumptions we use a dual form.

\begin{lemma}\label{lem:K}
Let $p,q$ be probability densities on $\R^d$, $g$ the density of $\mathcal N(b,\Sigma)$ with $\Sigma\succ0$, and $t\in[0,1]$. Put $k(y)=\int p(y-tv)\,q(y+(1-t)v)\,g(v)\,\dd v$. Then $\int\sqrt{k(y)}\,\dd y\le1$.
\end{lemma}

\begin{proof}
The substitution $v=u+b$, $y=y'+tb$ replaces $q$ by $q(\cdot+b)$ and $g$ by $g_\Sigma$, the density of $\mathcal N(0,\Sigma)$, so we may assume $b=0$. For $t=0$, $\int\sqrt k=\int\sqrt{p\,(q*g_\Sigma)}\le1$ by the Cauchy-Schwarz inequality, and $t=1$ is symmetric. Let $t\in(0,1)$ and let $g_{s\Sigma}$ be the density of $\mathcal N(0,s\Sigma)$. For $a,b',y\in\R^d$, expanding the quadratic forms shows
\begin{equation}\label{eq:bridge}
g_{t\Sigma}(y-a)\,g_{(1-t)\Sigma}(b'-y)=g_\Sigma(b'-a)\,g_{t(1-t)\Sigma}\big(y-(1-t)a-tb'\big).
\end{equation}
Indeed, with $D=b'-a$ and $y=(1-t)a+tb'+u$ one has $y-a=tD+u$ and $b'-y=(1-t)D-u$, and $\frac1t|tD+u|^2+\frac1{1-t}|(1-t)D-u|^2=|D|^2+\frac{|u|^2}{t(1-t)}$ in the norm $|\cdot|_\Sigma$; the normalizing constants agree. This is the Brownian bridge identity. Multiply \eqref{eq:bridge} by $p(a)q(b')$ and integrate over $(a,b')$. The map $(a,b')\mapsto((1-t)a+tb',b'-a)$ has unit Jacobian and $k$ is the density of the first coordinate under $p(a)q(b')g_\Sigma(b'-a)\,\dd a\,\dd b'$, so by Tonelli's theorem and the symmetry of $g_{(1-t)\Sigma}$,
\[
(k*g_{t(1-t)\Sigma})(y)=(p*g_{t\Sigma})(y)\,(q*g_{(1-t)\Sigma})(y).
\]
Since the square root is concave, $(\sqrt k)*g_{t(1-t)\Sigma}\le\sqrt{k*g_{t(1-t)\Sigma}}$ pointwise. Integrating and applying the Cauchy-Schwarz inequality,
\[
\int\sqrt k\le\int\sqrt{(p*g_{t\Sigma})(q*g_{(1-t)\Sigma})}\le1.
\]
\end{proof}

\begin{lemma}\label{lem:ac}
If $\mathcal E_g(0)$ and $\mathcal E_g(1)$ are finite, then $\Law(X)$, $\Law(X+V)$ and $\Law(X,V)$ have Lebesgue densities.
\end{lemma}

\begin{proof}
Let $\Pi=\Law(X,V)$, $\alpha=\Law(X)$, $\beta=\Law(X+V)$ and $S_1(x,v)=(x+v,v)$. Finiteness of $\mathcal E_g(0)$ gives $\Pi\ll\alpha\otimes g$, and finiteness of $\mathcal E_g(1)$ gives $\Pi\circ S_1^{-1}\ll\beta\otimes g$, that is, $\Pi\ll(\beta\otimes g)\circ S_1$. Let $\mathcal Z\subset\R^d$ be Lebesgue-null. Then $(\beta\otimes g)\circ S_1(\mathcal Z\times\R^d)=\int\!\!\int1_{\mathcal Z}(y-v)g(v)\,\dd v\,\beta(\dd y)=0$, so $\alpha(\mathcal Z)=0$. Hence $\alpha$ has a density, and so does $\beta$ by the symmetric argument. Then $\Pi\ll\alpha\otimes g$ is absolutely continuous.
\end{proof}

\begin{proof}[Second proof of Proposition~\ref{thm:entropy}, without moment assumptions]
Fix $t\in(0,1)$; the endpoints are trivial. By Lemma~\ref{lem:ac}, let $p$, $q$ and $\rho$ be the densities of $X$, $X+V$ and $(X,V)$. The pair $(M_t,V)$ has density $\rho_t(y,v)=\rho(y-tv,v)$, and $M_t$ has density $a(y)=\int\rho_t(y,v)\,\dd v$. All densities below are evaluated at points where they are positive and finite, which holds almost surely, and on $\{a=0\}$ the integrand in $v$ vanishes. Since a finite relative entropy has an integrable log-likelihood ratio,
\[
\mathcal E_g(0)=\E\log\frac{\rho(X,V)}{p(X)g(V)},\qquad\mathcal E_g(1)=\E\log\frac{\rho(X,V)}{q(X+V)g(V)},
\]
and put $\mathcal R=\dfrac{p(M_t-tV)\,q(M_t+(1-t)V)\,g(V)}{\rho_t(M_t,V)\,a(M_t)}$. By the Cauchy-Schwarz inequality in $v$ for fixed $y$, and Lemma~\ref{lem:K},
\[
\E\sqrt{\mathcal R}=\int\!\!\int\sqrt{\rho_t(y,v)}\sqrt{\frac{p(y-tv)q(y+(1-t)v)g(v)}{a(y)}}\,\dd v\,\dd y\le\int\sqrt{a}\sqrt{\frac{k}{a}}\,\dd y\le1 .
\]
Pointwise, the log-likelihood ratio $\log\frac{\rho_t(M_t,V)}{a(M_t)g(V)}$ defining $\mathcal E_g(t)$ is the sum of the two endpoint log-likelihood ratios and $\log\mathcal R$. The endpoint terms are integrable, and $(\log\mathcal R)^+\le2\sqrt{\mathcal R}$, so the positive part of this sum is integrable. Its negative part is integrable because it is a log-likelihood ratio between probability laws. Thus $\mathcal E_g(t)<\infty$, and taking expectations gives
\[
\mathcal E_g(t)=\mathcal E_g(0)+\mathcal E_g(1)+\E\log\mathcal R .
\]
Finally, $\log x\le2(\sqrt x-1)$ yields $\E\log\mathcal R\le2(\E\sqrt{\mathcal R}-1)\le0$, which proves the claim.
\end{proof}

\begin{remark}
The chain rule gives $\mathcal E_g(t)=\KL(\Law V\,\|\,g)+I(X+tV;V)$, with $I$ the mutual information. When the differential entropies involved are finite, choosing $g$ with the mean and covariance of $V$ gives the equivalent inequality
\[
\mathrm h(X+tV)\le I(X;X+V)+\mathrm h(\mathcal N_V),
\]
where $\mathcal N_V$ is that Gaussian. We are not aware of this endpoint inequality in the Brascamp-Lieb or sumset literature; see \cite{BW22,CL21,KM14,LCCV18} for related entropy inequalities. The Gaussian relative-entropy term cannot generally be omitted: $I(X+tV;V)$ can exceed $I(X;V)+I(X+V;V)$.\footnote{Let $V,W$ be independent and uniformly distributed on $\{-1,1\}$, and $X=VW$. Then $X$ is independent of $V$, so $I(X;V)=0$; $X+V=V(W+1)\in\{0,2V\}$ reveals $V$ with probability $\frac12$, so $I(X+V;V)=\frac12\log2$; and for $t\in(0,1)$, $|W+t|$ determines $W$, so $X+tV=V(W+t)$ determines $V$ and $I(X+tV;V)=\log2$. Here $\KL(\Law V\,\|\,g)=\infty$.}
\end{remark}

\subsection{The form used for OBABO}

For a law $p$ on $\R^{2d}$ with position marginal $\nu$, let $\Ev(p)=\KL(p\,\|\,\nu\otimes\varphi)$.

\begin{corollary}\label{cor:freeflight}
Let $p$ be a law on $\R^{2d}$, $h>0$ and $p_s=pA_s$. If $\Ev(p_0)$ and $\Ev(p_h)$ are finite, then $\Ev(p_s)\le\Ev(p_0)+\Ev(p_h)$ for every $s\in[0,h]$.
\end{corollary}

\begin{proof}
Let $(X,V)\sim p$ and $W=hV$. The bijection $(y,v)\mapsto(y,hv)$ maps $\Law(X+sV,V)$ to $\Law(X+(s/h)W,W)$ and $\Law(X+sV)\otimes\varphi$ to $\Law(X+(s/h)W)\otimes\mathcal N(0,h^2\Id)$, and relative entropy is invariant under bijections. Apply Proposition~\ref{thm:entropy} to the pair $(X,W)$ with $g$ the density of $\mathcal N(0,h^2\Id)$ and $t=s/h$.
\end{proof}

\section{Discrete entropic hypocoercivity for OBABO}\label{sec:hypo}

The Ornstein-Uhlenbeck half-steps dissipate velocity entropy. To obtain dissipation of position entropy, we add Lu's transport corrector and estimate its change under velocity Verlet. We retain the leading energy error as a signed difference, whose accumulated contribution is controlled by summation by parts. Assumption (A2) is used only in the proof of Theorem~\ref{thm:KL}, to bound the remainder $\tau$.

\subsection{Functionals}

Let $p$ be a law on $\R^{2d}$ with finite second moments and finite differential entropy, so that $\KL(p\|\mu)<\infty$, and let $\nu$ be its position marginal. Put
\[
\Ent(p)=\KL(p\|\mu),\quad\Ex(p)=\KL(\nu\|\mu_x),\quad\Ev(p)=\Ent(p)-\Ex(p)=\KL(p\|\nu\otimes\varphi),
\]
the last identity being the chain rule for relative entropy. We also write $\Ex(\nu)=\KL(\nu\|\mu_x)$ for a position law $\nu$. Since $\Ex(p)<\infty$, $\nu$ has a density, also denoted $\nu$, and its negative differential entropy $\Leb(\nu)=\int\nu\log\nu=-\mathrm h(\nu)$ is finite. We write $\Leb$ rather than $-\mathrm h$ to avoid a clash with the step size $h$. Let $T=T_\nu$ be the Brenier map \cite{Bre91} pushing $\nu$ to $\mu_x$, that is, the $\nu$-a.e.\ unique gradient of a convex function with $T_\#\nu=\mu_x$. Define
\begin{align*}
\Ccor(p)&=\E_p\langle v,x-T(x)\rangle, & J(p)&=\E_p\big|\E_p[v\,|\,x]\big|^2,\\
\cgap(\nu)&=\E_\nu\big[U(T(x))-U(x)-\langle\nabla U(x),T(x)-x\rangle\big]\ge0, & \Delta\Leb(\nu)&=\Leb(\nu)-\Leb(\mu_x).
\end{align*}
$\Ccor$ is Lu's transport corrector \cite[(2.18)]{Lu26}: along a free flight $p_s=pA_s$ with position marginal $\nu_s$, it is the derivative of $\frac12W_2^2(\nu_s,\mu_x)$ in $s$, for almost every $s$ along a sufficiently regular evolution \cite[(4.6)]{Lu26}.

\begin{lemma}[\cite{Lu26}]\label{lem:static}
Let $p$ be as above, with position marginal $\nu$. Then:
(S1) $J(p)\le2\Ev(p)$. (S2) $W_2^2(\nu,\mu_x)\le\frac2m\Ex(p)$. (S3) $|\Ccor(p)|\le\sqrt{J(p)}\,W_2(\nu,\mu_x)\le2\sqrt{\Ev(p)\Ex(p)/m}\le\Ent(p)/\sqrt m$.
\end{lemma}

\begin{proof}
(S1) is \cite[Lemma~3.1]{Lu26}, (S2) is Talagrand's inequality as in \cite[(3.10)]{Lu26}, and (S3) is \cite[Lemma~3.2]{Lu26}; we include the short proofs. (S1): for $\nu$-a.e.\ $x$, the Donsker-Varadhan formula with test functions $v\mapsto\langle a,v\rangle$ gives $\KL(p(\cdot|x)\,\|\,\varphi)\ge\sup_a(\langle a,\E[v|x]\rangle-|a|^2/2)=\frac12|\E[v|x]|^2$; integrate in $x$. (S2) is Talagrand's inequality for the $m$-strongly log-concave $\mu_x$ \cite{OV00}. (S3): condition on $x$ and apply the Cauchy-Schwarz inequality, using $\E_\nu|x-T(x)|^2=W_2^2$, then (S1), (S2) and $2\sqrt{ab}\le a+b$.
\end{proof}

\subsection{Exact one-step identities}

\begin{lemma}\label{lem:Ostep}
Let $p'=pO_{h/2}$. Then $p'$ and $p$ have the same position marginal, $\Ex(p')=\Ex(p)$, $\Ev(p')\le r^2\Ev(p)$ and $\Ccor(p')=r\,\Ccor(p)$.
\end{lemma}

\begin{proof}
$O_{h/2}$ acts on each fibre $\{x\}\times\R^d$ as the Ornstein-Uhlenbeck semigroup with invariant law $\varphi$ at time $\sqrt mh$, which contracts relative entropy to $\varphi$ by the factor $e^{-2\sqrt mh}=r^2$ (Bakry-\'Emery \cite{BGL14}). Integrating over $x$ gives the bound on $\Ev$. Since $\E[rv+\sigma\xi\,|\,x]=r\E[v|x]$ and $T$ depends only on $\nu$, $\Ccor(p')=r\Ccor(p)$.
\end{proof}

\begin{lemma}[Kicks]\label{lem:kick}
For $s\ge0$, $pB_s$ has the same position marginal $\nu$ as $p$, and
\[
\Ccor(pB_s)=\Ccor(p)-s\big[\cgap(\nu)+\Ex(\nu)-\Delta\Leb(\nu)\big].
\]
\end{lemma}

\begin{proof}
$\Ccor(pB_s)=\Ccor(p)-s\E_\nu\langle\nabla U,x-T(x)\rangle$. Pointwise, $\langle\nabla U(x),x-T(x)\rangle=c_T(x)+U(x)-U(T(x))$ with $c_T(x)=U(T(x))-U(x)-\langle\nabla U(x),T(x)-x\rangle$, and $T_\#\nu=\mu_x$, so $\E_\nu\langle\nabla U,x-T\rangle=\cgap(\nu)+\E_\nu U-\E_{\mu_x}U$. Finally $\Ex(\nu)=\Leb(\nu)+\E_\nu U+\log Z_x$ and $\Leb(\mu_x)=-\E_{\mu_x}U-\log Z_x$, where $Z_x=\int e^{-U}$, so $\E_\nu U-\E_{\mu_x}U=\Ex(\nu)-\Delta\Leb(\nu)$.
\end{proof}

\begin{lemma}[Volume-preserving maps]\label{lem:PhiEnt}
Let $\Phi'$ be a Borel bijection of $\R^{2d}$ with Borel inverse that preserves Lebesgue measure. If $p$ has finite second moments, $\Ent(p)<\infty$ and $\E_p|H\circ\Phi'-H|<\infty$, then $\Ent(p\Phi')=\Ent(p)+\E_p[H\circ\Phi'-H]$. This applies to $\Phi'=\Phi$.
\end{lemma}

\begin{proof}
Since $p$ has finite second moments, $\E_pH<\infty$, so $\int p\log p=\Ent(p)-\E_pH-\log\int e^{-H}$ is finite, and $\E_{p\Phi'}H=\E_pH\circ\Phi'$ is finite. $p\Phi'$ has density $p\circ\Phi'^{-1}$, so $\int p\Phi'\log p\Phi'=\int p\log p$. Hence $\Ent(p\Phi')-\Ent(p)=\E_p[H\circ\Phi'-H]$. $\Phi$ is a composition of symplectic maps, so it preserves Lebesgue measure.
\end{proof}
\subsection{Free transport}

Along a free flight $p_s=pA_s$, with $\nu_s$ the position marginal of $p_s$, two facts will be used, for every law $p$ with finite second moments and finite entropy. First, since $A_s$ preserves Lebesgue measure and $\E|v|^2$,
\begin{equation}\label{eq:EvLeb}
\Ev(p_s)=\int p\log p-\Leb(\nu_s)+\tfrac12\E_p|v|^2+\tfrac d2\log(2\pi),
\end{equation}
so $\Ev(p_s)+\Leb(\nu_s)$ does not depend on $s$,
because $\Ev(q)=\int q\log q-\int\nu_q\log\nu_q+\frac12\E_q|v|^2+\frac d2\log(2\pi)$ for every such law $q$ with position marginal $\nu_q$. Second, $\nu_s$ has a density for every $s$, because $X_s=X+sV$ has a density given $V$ when $(X,V)$ has one.

For a law $p$ with finite second moments and finite entropy write $p_s=pA_s$ and $\nu_s$ for the position marginal of $p_s$. By \eqref{eq:EvLeb}, $-\Delta\Leb(\nu_s)=\Ev(p_s)-\Ev(p)-\Delta\Leb(\nu_0)$, so for $t>0$
\begin{equation}\label{eq:freeint}
\int_0^t\big[6\Ev(p_s)-\Delta\Leb(\nu_s)\big]\dd s=7\int_0^t\Ev(p_s)\dd s-t\,\Ev(p)-t\,\Delta\Leb(\nu_0).
\end{equation}

\begin{lemma}[Semicontinuity]\label{lem:semicont}
Let $Y_k\to Y$ weakly in $\R^n$ with $\E|Y_k|^2\to\E|Y|^2<\infty$. Then $\limsup_k\mathrm h(Y_k)\le\mathrm h(Y)$. If $p_k\to p$ weakly on $\R^{2d}$, then $\liminf_k\Ev(p_k)\ge\Ev(p)$. Here differential entropy may equal $-\infty$, and we assign this value to laws without a Lebesgue density.
\end{lemma}

\begin{proof}
$\mathrm h(Y)=-\KL(\Law Y\|\mathcal N(0,\Id))+\frac12\E|Y|^2+\frac n2\log(2\pi)$, and relative entropy is jointly lower semicontinuous for weak convergence. For $\Ev$: if $p_k\to p$ weakly, then $\nu_k\otimes\varphi\to\nu\otimes\varphi$ weakly.
\end{proof}

\begin{lemma}[Stability of the corrector]\label{lem:corrstab}
Let $p_k\to p$ in $W_2(\R^{2d})$, and assume that the position marginal $\nu$ of $p$ has a density. For each $k$ let $\gamma_k$ be an optimal plan for $W_2$ between the position marginal $\nu_k$ of $p_k$ and $\mu_x$, and let $\Pi_k$ be the law of $(X,V,Y)$ with $(X,V)\sim p_k$, $(X,Y)\sim\gamma_k$ and $V,Y$ conditionally independent given $X$. Then $\int\langle v,x-y\rangle\,\dd\Pi_k\to\Ccor(p)$.
\end{lemma}

\begin{proof}
The family $(\Pi_k)$ is tight and has uniformly integrable second moments, because its marginals converge in $W_2$ or are fixed. Let $\Pi$ be a weak limit point. Its $(X,V)$ marginal is $p$, and its $(X,Y)$ marginal is a weak limit of optimal plans, hence optimal for $(\nu,\mu_x)$ \cite[Theorem~5.20]{Vil09OT}. Since $\nu$ has a density, this plan is unique and equals $(\mathrm{id},T)_\#\nu$ \cite{Bre91}. Hence $Y=T(X)$ $\Pi$-a.s., so $\Pi$ is the law of $(X,V,T(X))$ under $p$, and the whole sequence converges. The function $\langle v,x-y\rangle$ is continuous and bounded by $|v|^2+|x|^2+|y|^2$, which is uniformly integrable, so the integrals converge to $\Ccor(p)$.
\end{proof}

\begin{lemma}[Integrated free-transport inequality]\label{lem:approx}
Assume (A1). Let $p$ be a law on $\R^{2d}$ with finite second moments and finite entropy, $p_s=pA_s$, and $t>0$. Then
\begin{equation}\label{eq:freeinq}
\Ccor(p_t)-\Ccor(p)\le\int_0^t\big[6\Ev(p_s)-\Delta\Leb(\nu_s)\big]\dd s=7\int_0^t\Ev(p_s)\dd s-t\,\Ev(p)-t\,\Delta\Leb(\nu_0).
\end{equation}
\end{lemma}

\begin{proof}
We first establish the inequality for Gaussian convolutions of compactly supported laws, then approximate $p$ by such laws. Convergence in $W_2$ gives convergence of the correctors by Lemma~\ref{lem:corrstab}. For the right side, the representation \eqref{eq:EvLeb}, entropy semicontinuity, and entropy increase under Gaussian convolution give the upper bounds needed to pass to the limit.

\emph{The Gaussian-smoothed case.} Let $\tilde p$ be a law on $\R^{2d}$ with bounded support, $\iota>0$, $q=\tilde p*\mathcal N(0,\iota\Id_{2d})$, $q_s=qA_s$ and $\nu_s$ the position marginal of $q_s$. Then for every $t>0$,
\begin{equation}\label{eq:freesmooth}
\Ccor(q_t)-\Ccor(q)\le\int_0^t\big[6\Ev(q_s)-\Delta\Leb(\nu_s)\big]\dd s .
\end{equation}
This is an adaptation of Lu's bound on the time derivative of $\Ccor$ \cite[Lemma~4.1 and Section~5]{Lu26} to free transport. Lu proves that bound for solutions of the kinetic Langevin equation in the class of \cite[Definition~2.5]{Lu26}, which does not contain free transport, so we prove \eqref{eq:freesmooth} in Appendix~\ref{app:free}. The proof follows Lu's method: position marginals at nearby times are coupled through the flow of the conditional mean velocity, which bounds the second derivative of $\frac12W_2^2(\nu_s,\mu_x)$, and the term produced by the conditional velocity covariance is bounded by an estimate of Lu \cite[Lemmas~5.1 to~5.3]{Lu26}, a static inequality for Brenier maps, whose proof we include (Lemma~\ref{lem:stress}); the term $J(q_s)$ that this produces is bounded by (S1).

\emph{Approximants.} Let $B_R\subset\R^{2d}$ be the closed ball of radius $R$ about $0$. For $R$ with $p(B_R)>0$ and $\iota\in(0,1]$ let $p_R=p(\cdot\,|\,B_R)$ and $q_{R,\iota}=p_R*\mathcal N(0,\iota\Id_{2d})$, and write $\nu(q)$ for the position marginal of a law $q$. By \eqref{eq:freesmooth} with $\tilde p=p_R$ and \eqref{eq:freeint},
\[
\Ccor(q_{R,\iota}A_t)-\Ccor(q_{R,\iota})\le7\int_0^t\Ev(q_{R,\iota}A_s)\dd s-t\,\Ev(q_{R,\iota})-t\,\Delta\Leb(\nu(q_{R,\iota})).
\]
Let $R_k\to\infty$ and $\iota_k\to0$, and put $q_k=q_{R_k,\iota_k}$.

\emph{Left side.} $p_R\to p$ in total variation and $\int_{B_R}|z|^2p(\dd z)\to\int|z|^2p(\dd z)$, so $W_2(p_R,p)\to0$; also $W_2^2(q_{R,\iota},p_R)\le2d\iota$. Hence $q_k\to p$ in $W_2$, and since $A_t$ is linear, $q_kA_t\to p_t$ in $W_2$. The position marginals of $p$ and $p_t$ have densities, so Lemma~\ref{lem:corrstab} gives convergence of both correctors.

\emph{The integral.} Let $(X^k,V^k)\sim p_{R_k}$, $(X,V)\sim p$, and let $\zeta'=(\zeta_x,\zeta_v)\sim\mathcal N(0,\Id_{2d})$ be independent of $(X^k,V^k)$. By \eqref{eq:EvLeb}, and since entropy increases under independent addition,
\[
\Ev(q_kA_s)\le-\mathrm h(p_{R_k})+\mathrm h\big(X^k+sV^k+\sqrt{\iota_k}(\zeta_x+s\zeta_v)\big)+\tfrac12\E|V^k|^2+\tfrac{\iota_kd}2+\tfrac d2\log(2\pi).
\]
Here $\mathrm h(p_R)=-p(B_R)^{-1}\int_{B_R}p\log p+\log p(B_R)\to\mathrm h(p)$ by dominated convergence, since $p\log p$ is integrable (its negative part because $p$ has finite second moments), and $\E|V^k|^2\to\E|V|^2$. For each $s$, $X^k+sV^k+\sqrt{\iota_k}(\zeta_x+s\zeta_v)\to X+sV$ in $W_2$, so by Lemma~\ref{lem:semicont} and \eqref{eq:EvLeb} the $\limsup$ of the right side as $k\to\infty$ is at most $\Ev(p_s)$. The right side is bounded above uniformly in $(s,k)$, since $\mathrm h(Y)\le\frac d2\log(2\pi\mathrm e\,\E|Y|^2/d)$. By the reverse Fatou lemma, $\limsup_k\int_0^t\Ev(q_kA_s)\dd s\le\int_0^t\Ev(p_s)\dd s$.

\emph{Endpoint terms.} $\liminf_k\Ev(q_k)\ge\Ev(p)$ by Lemma~\ref{lem:semicont}. Next, $-\Delta\Leb(\nu(q_k))=\mathrm h(\nu(q_k))-\mathrm h(\mu_x)$, and $\limsup_k\mathrm h(\nu(q_k))\le\mathrm h(\nu_0)$ by Lemma~\ref{lem:semicont}.

\emph{Conclusion.} The left side converges to $\Ccor(p_t)-\Ccor(p)$, and the $\limsup$ of the right side is at most the right side of \eqref{eq:freeinq}. The equality in \eqref{eq:freeinq} is \eqref{eq:freeint}.
\end{proof}

No Hessian-Lipschitz condition is used in Lemma~\ref{lem:approx}. In continuous time the term $-\Delta\Leb$ from free transport cancels the term $+\Delta\Leb$ produced by the force (Lemma~\ref{lem:kick} with $s=\dd t$), since both act at the same instant. For OBABO the force acts at the two ends of the flight, and the cancellation becomes a trapezoidal rule.

\subsection{The trapezoidal rule error}

For one OBABO step from a law $p$ write
\[
p'=pO_{h/2},\qquad p_0=p'B_{h/2},\qquad p_s=p_0A_s\ (s\in[0,h]),\qquad p'\Phi=p_hB_{h/2},
\]
let $\nu_s$ be the position marginal of $p_s$, so that $\nu_0$ is also that of $p'$ and $\nu_h$ that of $p'\Phi$, and let $V$ denote the velocity under $p_0$, which equals $V(w)$ in \eqref{eq:halfvel} for $w\sim p'$. Define
\[
\Etrap=\tfrac h2\big[\Leb(\nu_0)+\Leb(\nu_h)\big]-\int_0^h\Leb(\nu_s)\,\dd s .
\]

\begin{lemma}\label{lem:identities}
Assume $p'$ has finite second moments and finite entropy.
\begin{enumerate}
\item $\Etrap=\int_0^h\Ev(p_s)\dd s-\frac h2[\Ev(p_0)+\Ev(p_h)]$.
\item With $\varrho_0=\Ev(p_0)-\Ev(p')$ and $\varrho_1=\Ev(p'\Phi)-\Ev(p_h)$,
\[
\varrho_0-\varrho_1=\tfrac{h^2}{2}\E\langle V,\bar AV\rangle-\E\big[G(x_0)+G(x_h)\big]\le\tfrac{\bar h^2}{2}\E|V|^2,
\]
where $\bar A=\int_0^1\nabla^2U(x_0+\vartheta hV)\dd\vartheta$ and $(x_0,V)\sim p_0$, $x_h=x_0+hV$.
\end{enumerate}
\end{lemma}

\begin{proof}
(i) Integrate $\Leb(\nu_s)=\mathrm{const}-\Ev(p_s)$, from \eqref{eq:EvLeb}.
(ii) A kick moves no position and is a shear, so it changes $\Ev$ by the change of $\E|v|^2/2$. With $v=V+\frac h2\nabla U(x_0)$ before the first kick and $V-\frac h2\nabla U(x_h)$ after the second,
\[
\varrho_0=-\tfrac h2\E\langle V,\nabla U(x_0)\rangle-\E G(x_0),\qquad\varrho_1=-\tfrac h2\E\langle V,\nabla U(x_h)\rangle+\E G(x_h).
\]
Subtract and use $\nabla U(x_h)-\nabla U(x_0)=h\bar AV$ with $0\preceq\bar A\preceq L\Id$.
\end{proof}

\begin{figure}[t]
\centering
\includegraphics[width=\textwidth]{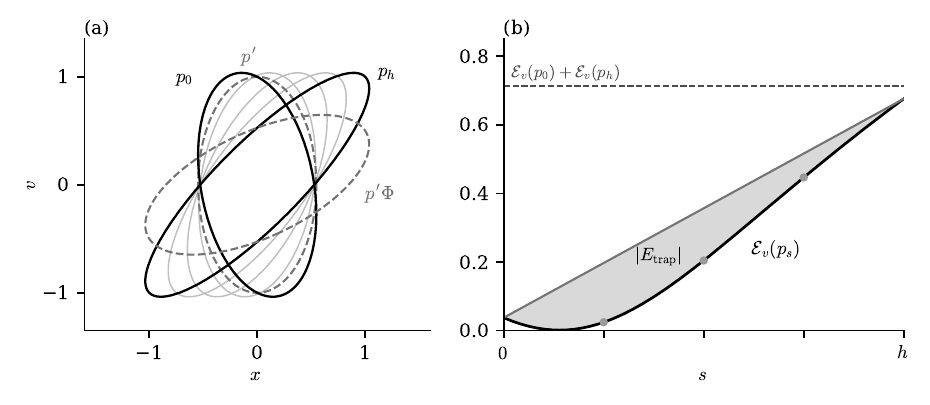}
\caption{One velocity Verlet step. (a) One-standard-deviation ellipses of the laws in the step; light gray: $p_s$ for $s=\frac h4,\frac h2,\frac{3h}4$. (b) $\Ev(p_s)$ (black), its trapezoidal chord (gray), and the bound of Corollary~\ref{cor:freeflight} (dashed); the shaded area is $|\Etrap|$.}
\label{fig:trap}
\end{figure}

By (i), $\Etrap$ is the error of the trapezoidal rule for $\int_0^h\Ev(p_s)\dd s$, and Corollary~\ref{cor:freeflight} bounds the integrand by $\Ev(p_0)+\Ev(p_h)$. Figure~\ref{fig:trap} illustrates the velocity Verlet stage $p'\mapsto p'\Phi$ between the two Ornstein-Uhlenbeck half-steps. We take $U(x)=\frac12x^2$ on $\R$, $h=1$ and $p'=\mathcal N(0,\mathrm{diag}(0.3,1))$; all laws shown are Gaussian, and $\Ev(p_s)$ is explicit.

\subsection{The corrector over one velocity Verlet step}

\begin{proposition}\label{prop:DS0}
If $p'$ has finite second moments and finite entropy, then
\begin{align*}
\Ccor(p'\Phi)-\Ccor(p')\le{}&6\int_0^h\Ev(p_s)\,\dd s-\tfrac h2\big[\Ex(\nu_0)+\Ex(\nu_h)\big]\\
&-\tfrac h2\big[\cgap(\nu_0)+\cgap(\nu_h)\big]+\Etrap .
\end{align*}
\end{proposition}

\begin{proof}
By Lemma~\ref{lem:kick}, $\Ccor(p_0)=\Ccor(p')-\frac h2[\cgap(\nu_0)+\Ex(\nu_0)-\Delta\Leb(\nu_0)]$ and $\Ccor(p'\Phi)=\Ccor(p_h)-\frac h2[\cgap(\nu_h)+\Ex(\nu_h)-\Delta\Leb(\nu_h)]$. By Lemma~\ref{lem:approx} applied to $p_0$, $\Ccor(p_h)-\Ccor(p_0)\le\int_0^h[6\Ev(p_s)-\Delta\Leb(\nu_s)]\dd s$. Add the three relations; the constants $\Leb(\mu_x)$ cancel and $\frac h2[\Delta\Leb(\nu_0)+\Delta\Leb(\nu_h)]-\int_0^h\Delta\Leb(\nu_s)\dd s=\Etrap$.
\end{proof}

\begin{proposition}[Corrector bound for one step]\label{prop:DS}
If $p'$ has finite second moments and finite entropy, then
\begin{equation}\label{eq:DS}
\Ccor(p'\Phi)-\Ccor(p')\le h\big[a_1\big(\Ev(p')+\Ev(p'\Phi)\big)-a_2\Ex(p')\big]+h\,b(p'),
\end{equation}
with $a_1=\frac{13}2$, $a_2=\frac12$ and $b(p')=\frac{13}{4}\bar h^2\E|V|^2$.
\end{proposition}

\begin{proof}
In Proposition~\ref{prop:DS0} drop the gaps and insert Lemma~\ref{lem:identities}(i):
\begin{equation}\label{eq:star}
\Ccor(p'\Phi)-\Ccor(p')\le7\int_0^h\Ev(p_s)\dd s-\tfrac h2\big[\Ev(p_0)+\Ev(p_h)\big]-\tfrac h2\big[\Ex(\nu_0)+\Ex(\nu_h)\big].
\end{equation}
By Corollary~\ref{cor:freeflight}, $\int_0^h\Ev(p_s)\dd s\le h[\Ev(p_0)+\Ev(p_h)]$, so the right side is at most $\frac{13}2h[\Ev(p_0)+\Ev(p_h)]-\frac h2\Ex(\nu_0)$. By Lemma~\ref{lem:identities}(ii), $\Ev(p_0)+\Ev(p_h)=\Ev(p')+\Ev(p'\Phi)+\varrho_0-\varrho_1\le\Ev(p')+\Ev(p'\Phi)+\frac{\bar h^2}2\E|V|^2$. Finally $\Ex(\nu_0)=\Ex(p')$.
\end{proof}

No Fisher information, no Hessian-Lipschitz constant and no lower bound on the time elapsed enter \eqref{eq:DS}.

\subsection{One-step decay of $\Ent+\eta\,\Ccor$}

The next proposition is stated for a class of maps that contains $\Phi$ and the map $\Phi_\sharp$ of Section~\ref{sec:comparison}.

\begin{proposition}\label{prop:R}
Let $\Phi'$ be a Borel bijection of $\R^{2d}$ with Borel inverse that preserves Lebesgue measure, such that $H(\Phi'(w))-H(w)=G(x'_1)-G(x)+\tau'(w)$ for all $w=(x,v)$, where $x'_1$ is the position of $\Phi'(w)$. Let $K'(z,\cdot)=\Law\big((O_{h/2}^{\xi'}\circ\Phi'\circ O_{h/2}^{\xi})(z)\big)$, where $\xi$ and $\xi'$ are independent random vectors, each with law $\varphi$, and let $p_k=p_0K'^k$, $p_k'=p_kO_{h/2}$, $g_k=\E_{p_k}G$ and $e_k=\E_{p_k'}[\tau']$.

Assume $\Ent(p_0)<\infty$, that each $p_k$ has finite third moments, and that for every $k$
\[
\Ccor(p_k'\Phi')-\Ccor(p_k')\le h\big[a_1\big(\Ev(p_k')+\Ev(p_k'\Phi')\big)-a_2\Ex(p_k')+b_k\big],
\]
where $a_1,a_2>0$ and $b_k\ge0$.

Let
\[
0<\theta\le\min\Big\{\frac{1-\sqrt mh}{2a_1},\ \frac{ra_2}{8},\ \frac{1-\sqrt mh}{ra_2},\ \frac12\Big\},\qquad\eta=\theta\sqrt m,\qquad\lambda=\frac{\theta ra_2\sqrt m}{2(1+\theta)},
\]
and $\mathcal H=\Ent+\eta\,\Ccor$.

Then $\mathcal H(p_{k+1})\le(1-\lambda h)\mathcal H(p_k)+(g_{k+1}-g_k)+e_k+\eta hb_k$ for every $k$.

Moreover, for $n\ge1$,
\[
\Ent(p_n)\le\frac{1+\theta}{1-\theta}e^{-\lambda nh}\Ent(p_0)+\frac1{1-\theta}\Big[g_n+\frac1{\lambda h}\sup_{k<n}\big(e_k^++\eta hb_k\big)\Big].
\]
\end{proposition}

\begin{proof}
Fix $k$ and write $p=p_k$, $p'=pO_{h/2}$, $p''=p'\Phi'$, so $p_{k+1}=p''O_{h/2}$, and $\delta=1-r^2=1-e^{-\gamma h}$. All quantities are finite: $\Ent(p_{k+1})\le\Ent(p_k)+\E_{p_k'}|H\circ\Phi'-H|$ by Lemmas~\ref{lem:Ostep} and~\ref{lem:PhiEnt}, and $H$ and $H\circ\Phi'$ are integrable under $p_k'$, since $p_k'$ and $p_k'\Phi'$ have finite second moments; hence so are $|H\circ\Phi'-H|$ and $\tau'$. Indeed, the first Ornstein-Uhlenbeck half-step gives $p_k'$ finite second moments. For the final half-step, writing the velocity of $p_{k+1}$ as $Y=rV+\sigma\xi'$ with $V$ the velocity of $p_k'\Phi'$ gives $r^2\E|V|^2\le2\E|Y|^2+2\sigma^2d$, and positions are unchanged, so $p_k'\Phi'$ also has finite second moments. Inductively, Lemmas~\ref{lem:Ostep} and~\ref{lem:PhiEnt} imply that all entropy and energy terms below are finite.

\emph{Entropy.} By Lemma~\ref{lem:Ostep}, $\Ent(p)-\Ent(p')=\Ev(p)-\Ev(p')\ge\delta\Ev(p)$ and $\Ent(p'')-\Ent(p_{k+1})\ge\delta\Ev(p'')$. By Lemma~\ref{lem:PhiEnt}, $\Ent(p'')=\Ent(p')+\E_{p'}[G(x_1')-G(x)]+e_k$, and $\E_{p'}[G(x_1')-G(x)]=g_{k+1}-g_k$ because the Ornstein-Uhlenbeck half-steps do not move positions. Hence
\[
\Ent(p_{k+1})\le\Ent(p)-\delta\big[\Ev(p)+\Ev(p'')\big]+(g_{k+1}-g_k)+e_k .
\]
\emph{Corrector.} By Lemma~\ref{lem:Ostep}, $\Ccor(p')=r\Ccor(p)$ and $\Ccor(p_{k+1})=r\Ccor(p'')$. With the assumed corrector bound, $\Ev(p')\le\Ev(p)$, $\Ex(p')=\Ex(p)$ and $r\le1$,
\[
\Ccor(p_{k+1})\le r^2\Ccor(p)+ha_1\big[\Ev(p)+\Ev(p'')\big]-rha_2\Ex(p)+hb_k .
\]
\emph{Combination.} By (S3), $-\delta\Ccor(p)\le2\delta\sqrt{\Ev(p)\Ex(p)/m}$, and $2\eta\delta\sqrt{\Ev\Ex/m}\le\frac\delta2\Ev+\frac{2\eta^2\delta}{m}\Ex$. Therefore
\begin{align*}
\mathcal H(p_{k+1})-\mathcal H(p)\le{}&\Big(-\tfrac\delta2+\eta ha_1\Big)\Ev(p)+\big(-\delta+\eta ha_1\big)\Ev(p'')+\eta\Big(\tfrac{2\theta\delta}{\sqrt m}-rha_2\Big)\Ex(p)\\
&+(g_{k+1}-g_k)+e_k+\eta hb_k .
\end{align*}
Since $\gamma h(1-\gamma h/2)\le\delta\le\gamma h$, $\gamma=2\sqrt m$ and $\gamma h/2=\sqrt mh$, the first constraint on $\theta$ gives $\eta ha_1\le\delta/4$, and the second gives $2\theta\delta/\sqrt m\le4\theta h\le rha_2/2$. So the coefficient of $\Ev(p)$ is at most $-\delta/4$, that of $\Ev(p'')$ is nonpositive, and that of $\Ex(p)$ is at most $-\theta ra_2h\sqrt m/2$. By the third constraint, $\delta/4\ge\frac12(1-\sqrt mh)h\sqrt m\ge\theta ra_2h\sqrt m/2$, so $\mathcal H(p_{k+1})-\mathcal H(p)\le-\frac{\theta ra_2h\sqrt m}{2}\Ent(p)+(g_{k+1}-g_k)+e_k+\eta hb_k$. By (S3), $(1-\theta)\Ent\le\mathcal H\le(1+\theta)\Ent$, which gives the one-step inequality.

\emph{Unrolling.} The parameter restrictions imply $0<\lambda h<1$, so the geometric weights below are nonnegative. With $w_k=(1-\lambda h)^{n-1-k}$, iteration gives $\mathcal H(p_n)\le(1-\lambda h)^n\mathcal H(p_0)+\sum_{k<n}w_k[(g_{k+1}-g_k)+e_k+\eta hb_k]$. Summation by parts gives $\sum_{k<n}w_k(g_{k+1}-g_k)=g_n-w_0g_0-\sum_{k=1}^{n-1}(w_k-w_{k-1})g_k\le g_n$, because $g_k\ge0$ and $w_k$ increases with $w_{n-1}=1$. Use $\sum_kw_k\le(\lambda h)^{-1}$, $(1-\lambda h)^n\le e^{-\lambda nh}$ and the comparison of $\mathcal H$ with $\Ent$.
\end{proof}

\begin{remark}
The leading energy error contributes the differences $g_{k+1}-g_k$. Summation by parts leaves $g_n$ and nonpositive terms, so their geometrically weighted contribution is at most $g_n\le\frac{\bar h^2}4\E_{p_n}(U-U_*)$. This is the same telescoping that the non-reversibly updated uniform exploits (Proposition~\ref{prop:allaccept}).
\end{remark}

\subsection{Proof of Theorem~\ref{thm:KL}}

\begin{lemma}[Starting entropy]\label{lem:start}
Under (A1), $\KL(\delta_zK\|\mu)\le\E_zH(Z_1)+d\log\dfrac{\kappa}{2e\bar h^2(1-\sqrt mh)}$ for $\sqrt mh<1$, and the density of $\delta_zK$ is bounded by $(2\pi h\sigma^2)^{-d}$.
\end{lemma}

\begin{proof}
From $z=(x,v)$, $X_1=x+hrv-\frac{h^2}2\nabla U(x)+h\sigma\xi$, and given $X_1$ the velocity before the second Ornstein-Uhlenbeck half-step is determined, so that $Z_1$ has the density of $X_1\sim\mathcal N(x+hrv-\frac{h^2}{2}\nabla U(x),h^2\sigma^2\Id)$ and $V_1|X_1\sim\mathcal N(r(\frac{X_1-x}h-\frac h2\nabla U(X_1)),\sigma^2\Id)$. This density is bounded by $(2\pi h^2\sigma^2)^{-d/2}(2\pi\sigma^2)^{-d/2}$, and its differential entropy is $\frac d2\log(2\pi eh^2\sigma^2)+\frac d2\log(2\pi e\sigma^2)$. With $\int e^{-H}\le(2\pi)^{d}m^{-d/2}$,
\[
\KL(\delta_zK\|\mu)=\E_zH(Z_1)+\log\int e^{-H}-\mathrm h(Z_1)\le\E_zH(Z_1)+\tfrac d2\log\frac1{e^2mh^2\sigma^4}.
\]
Finally $\sigma^2=1-e^{-2\sqrt mh}\ge2\sqrt mh(1-\sqrt mh)$ and $mh^2=\bar h^2/\kappa$.
\end{proof}

\begin{proof}[Proof of Theorem~\ref{thm:KL}]
Apply Proposition~\ref{prop:R} with $\Phi'=\Phi$ to $p_0=\delta_zK$, so that $p_k=\Law(Z_{k+1})$. The laws $p_k$ have finite moments of all orders, since $\nabla U$ is Lipschitz, so that each of the maps $B_s$, $A_s$ has at most linear growth, and the noise is Gaussian, and finite entropy, since $\delta_zK$ has a bounded density and $K$ does not increase entropy relative to $\mu$ by more than $\E|\Delta H|$. By Proposition~\ref{prop:DS}, \eqref{eq:DS} holds at $p_k'=\Law(W_{k+1})$ with $a_1=\frac{13}2$, $a_2=\frac12$, $b_k=\frac{13}4\bar h^2\E|V_{k+1}|^2$. Take $\theta=r/16$. The second constraint holds with equality, the first holds because $\sqrt mh\le3/16$ gives $(1-\sqrt mh)/13\ge1/16\ge\theta$, and the third because $2(1-\sqrt mh)/r\ge1$. Then $\frac{1+\theta}{1-\theta}\le\frac{17}{15}$, $\frac1{1-\theta}\le\frac{16}{15}$, $\lambda$ is \eqref{eq:lambda}, and $\eta/\lambda=(16+r)/(4r)$, so $\frac{\eta}{\lambda}\cdot\frac{13}{4}\le\frac{221}{16r}$. By Lemma~\ref{lem:energy}, $e_k^+\le\frac{L_H}{12}h^3\E|V_{k+1}|^3$, and $g_n\le\frac{\bar h^2}4\E(U(X_{n+1})-U_*)$. Lemma~\ref{lem:start} bounds $\Ent(p_0)$. The lower bound in \eqref{eq:lambda} follows from $r\ge e^{-\sqrt mh}$ and $16+r\le17$.
\end{proof}

\section{The non-reversibly updated uniform, the coupling argument and moment bounds}\label{sec:coupling}

Proposition~\ref{prop:allaccept} bounds rejection over an epoch using the accumulated signed energy error. Proposition~\ref{prop:coupling} then combines an endpoint coupling of the unadjusted chains with this rejection bound. Section~\ref{sec:moments} controls the unadjusted chain while friction dominates the Verlet error. Section~\ref{sec:comparison} introduces a comparison chain to handle faster flights and apply the entropy estimate of Section~\ref{sec:hypo}.

\subsection{Acceptance of every move of an epoch}

Consider one epoch of Definition~\ref{alg:pu} from $z$, driven by noises $(\xi_k,\xi_k')_{k<E}$ and the exponential variable $\ell$. Let $(Z_k,W_k)$ be the unadjusted chain from $z$ driven by the same noises, and put
\[
S_n=\sum_{k<n}\Delta H(W_k)=G(X_n)-G(X_0)+\sum_{k<n}\tau(W_k),\qquad 1\le n\le E,
\]
the second form by Lemma~\ref{lem:energy}, since the Ornstein-Uhlenbeck half-steps do not move positions. After $n$ accepted moves, the remaining exponential threshold is $\ell-S_n$. Thus acceptance throughout the epoch depends on the maximum of the partial sums, rather than on the sum of the positive one-step errors.

\begin{proposition}\label{prop:allaccept}
Every move of the epoch is accepted if and only if $\ell>\max_{1\le n\le E}S_n$, and on this event the adjusted epoch coincides with the unadjusted chain. Consequently
\begin{align*}
\PP(\text{a move of the epoch is rejected})&=\E\big[1-e^{-\max_nS_n^+}\big]\\
&\le\E\Big[\min\Big(1,\tfrac{\bar h^2}4\max_{n\le E}\big(U(X_n)-U_*\big)+\sum_{k<E}|\tau(W_k)|\Big)\Big].
\end{align*}
\end{proposition}

\begin{proof}
By induction, if the first $n$ moves are accepted, the adjusted path equals the unadjusted one up to step $n$ and $\ell$ has become $\ell-S_n$; the next move is accepted if and only if $\ell-S_n>\Delta H(W_n)$, that is, $\ell>S_{n+1}$. Since $\ell$ is independent of the noises, the conditional probability of a rejection given the noises is $1-e^{-\max_nS_n^+}\le\min(1,\max_nS_n^+)$. Finally $S_n\le G(X_n)+\sum_{k<n}|\tau(W_k)|$ and $G\le\frac{\bar h^2}4(U-U_*)$.
\end{proof}

The leading energy error contributes a maximum over visited positions, while the remainder $\sum_{k<E}|\tau(W_k)|$ still accumulates over the epoch. With a fresh uniform at every step the corresponding bound is $\sum_{k<E}\E[1-e^{-\Delta H(W_k)^+}]$, a sum of $E$ terms. For quadratic $U$, $\tau\equiv0$ and $\ell+G(x)$ is invariant during an epoch, so the adjusted epoch is unadjusted OBABO with the velocity reversed whenever a move would leave $\{G<G(x_0)+\ell\}$.

\subsection{The coupling theorem}

For $z\in\R^{2d}$ let $\beta(z)$ denote the probability in Proposition~\ref{prop:allaccept}. The next lemma couples the noises of a whole epoch rather than step by step. This is necessary: for Gaussian targets, the meeting probability of any Markovian coupling of two copies of OBABO decays strictly slower than their total variation distance \cite{BCS26}.

\begin{lemma}[Maximal endpoint coupling]\label{lem:glue}
There is a probability kernel $(z,\tilde z)\mapsto\Gamma_{z,\tilde z}$ from $\R^{2d}\times\R^{2d}$ to pairs of noise blocks $((\xi_k,\xi_k')_{k<E},(\tilde\xi_k,\tilde\xi_k')_{k<E})$ such that both marginals are i.i.d.\ standard Gaussian and, for the unadjusted chains $Z_E$ from $z$ and $\tilde Z_E$ from $\tilde z$ driven by the two blocks, $\Gamma_{z,\tilde z}(Z_E\ne\tilde Z_E)=\mathrm{TV}(\delta_zK^E,\delta_{\tilde z}K^E)$.
\end{lemma}

\begin{proof}
We first maximally couple the endpoint distributions. Conditional on those endpoints, we sample the two paths from their respective conditional path laws. We then recover each Gaussian noise block from its path. By the proof of Lemma~\ref{lem:start}, $K(z,\dd z')=k(z,z')\dd z'$ with $k$ continuous and strictly positive on $\R^{2d}\times\R^{2d}$. The path $(Z_1,\dots,Z_E)$ from $z$ has density $q_z(y_{1:E})=\prod_{i=1}^Ek(y_{i-1},y_i)$ with $y_0=z$, the endpoint has density $f_z(y)=\int q_z(y_{1:E-1},y)\,\dd y_{1:E-1}$, and $\beta_{z,y}(\dd y_{1:E-1})=q_z(y_{1:E-1},y)\dd y_{1:E-1}/f_z(y)$ is the law of the path given $Z_E=y$; all are jointly measurable by Tonelli's theorem. Let $a=\int f_z\wedge f_{\tilde z}=1-\mathrm{TV}(\delta_zK^E,\delta_{\tilde z}K^E)$. On pairs of paths define
\[
\int(f_z\wedge f_{\tilde z})(y)\,\beta_{z,y}\otimes\beta_{\tilde z,y}\,\dd y+\frac1{1-a}\iint(f_z-f_z\wedge f_{\tilde z})(y)(f_{\tilde z}-f_z\wedge f_{\tilde z})(\tilde y)\,\beta_{z,y}\otimes\beta_{\tilde z,\tilde y}\,\dd y\,\dd\tilde y,
\]
with paths ending at $y$ and $\tilde y$ respectively, and without the second term if $a=1$. This is a probability kernel in $(z,\tilde z)$ whose marginals are the path laws from $z$ and $\tilde z$, and the endpoints differ with probability $1-a$. Given $z$, the map from noises to paths, $(\xi_k,\xi_k')_{k<E}\mapsto(Z_k)_{k\le E}$, is a bimeasurable bijection: $\xi_k=(v'_k-rv_k)/\sigma$ with $v'_k=(x_{k+1}-x_k)/h+\frac h2\nabla U(x_k)$, and $\xi'_k=(v_{k+1}-r(v'_k-\frac h2\nabla U(x_k)-\frac h2\nabla U(x_{k+1})))/\sigma$. Transport the kernel through it.
\end{proof}

\begin{proposition}[Coupling]\label{prop:coupling}
Let $\mathcal D\subset\R^{2d}$ be measurable and
\[
q=\sup_{z,\tilde z\in \mathcal D}\Big(\mathrm{TV}(\delta_zK^E,\delta_{\tilde z}K^E)+\beta(z)+\beta(\tilde z)\Big).
\]
If $M_2=\|\dd\nu/\dd\mu\|_{L^2(\mu)}<\infty$, then for every $N\ge1$, $\mathrm{TV}(\nu P^N,\mu)\le q^N+N\mu(\mathcal D^c)+M_2\big(N\mu(\mathcal D^c)\big)^{1/2}$.
\end{proposition}

\begin{proof}
Run two copies $(\mathcal Z_t)$ from $\nu$ and $(\tilde{\mathcal Z}_t)$ from $\mu$, epoch by epoch. At the start of an epoch with states $(z,\tilde z)$, both copies use the same exponential variable $\ell$, drawn independently of everything else. If $z=\tilde z$ or $(z,\tilde z)\notin \mathcal D\times \mathcal D$ they use the same noises. Otherwise they use noise blocks drawn from $\Gamma_{z,\tilde z}$ of Lemma~\ref{lem:glue}, independently of the past. Given the past, each copy's noises are i.i.d.\ standard Gaussian and independent of $\ell$, so each copy is a version of the adjusted chain, and once the copies meet at an epoch boundary they stay equal.

For $z,\tilde z\in \mathcal D$ with $z\ne\tilde z$, on the event that neither copy rejects a move, both copies coincide with their unadjusted chains (Proposition~\ref{prop:allaccept}). The probability that the first copy rejects some move depends only on its own law and equals $\beta(z)$, and likewise for the second. Hence the probability that the copies differ at the end of the epoch is at most $q$.

Let $\mathcal B$ be the event that $\mathcal Z_{jE}\notin \mathcal D$ or $\tilde{\mathcal Z}_{jE}\notin \mathcal D$ for some $j<N$. Under $\PP_\mu$ each $\mathcal Z_{jE}$ has law $\mu$ (Lemma~\ref{lem:invariance}), so $\PP_\mu(\exists j<N:\mathcal Z_{jE}\notin \mathcal D)\le N\mu(\mathcal D^c)$, and the same bound holds for the second copy. For the first copy, started from $\nu$, the Cauchy-Schwarz inequality gives $\PP_\nu(\mathcal A)=\E_\mu[\frac{\dd\nu}{\dd\mu}(\mathcal Z_0)\PP(\mathcal A|\mathcal Z_0)]\le M_2\PP_\mu(\mathcal A)^{1/2}$ for every event $\mathcal A$ of its path. So $\PP(\mathcal B)\le N\mu(\mathcal D^c)+M_2(N\mu(\mathcal D^c))^{1/2}$. On $\mathcal B^c\cap\{\mathcal Z_{NE}\ne\tilde{\mathcal Z}_{NE}\}$ the copies differ at every epoch boundary $jE\le NE$ and start every epoch in $\mathcal D\times \mathcal D$. The pair of copies is a Markov chain on $\R^{4d}$ with a measurable epoch kernel, since $\mathcal D$ and the diagonal are measurable and $\Gamma$ is a kernel. Let $A_j$ be the event that for every $i<j$ the copies are distinct and in $\mathcal D\times \mathcal D$ at time $iE$ and distinct at time $(i+1)E$. The previous paragraph gives $\PP(A_{j+1}\,|\,\mathcal F_{jE})\le q1_{A_j}$, with $\mathcal F_{jE}$ the $\sigma$-field generated by both paths up to time $jE$, hence $\PP(A_N)\le q^N$, and the event above is contained in $A_N$. The coupling inequality gives the claim.
\end{proof}

Here localization is needed only at epoch boundaries, rather than at every step as in \cite[Assumption~A.S(iv)]{BRO24}; stationarity and the $L^2(\mu)$ warm start give the required bound directly. The coupling of Lemma~\ref{lem:glue} is not adapted: the noises of one copy at a given step depend on the entire noise block of the other copy. This is harmless because the rejection probabilities $\beta(z)$ and $\beta(\tilde z)$ depend only on the marginal law of each copy's noises.

\begin{lemma}[High-energy mass]\label{lem:tail}
Under (A1), for every $R\ge0$, $\mu(H>R)\le2^{3d/2}e^{-R/2}$.
\end{lemma}

\begin{proof}
Assume $x_*=0$ and $U_*=0$. For $t\in(0,1]$, convexity gives $U(y)\le tU(y/t)$, so $\int e^{-tU(x)}\dd x=t^{-d}\int e^{-tU(y/t)}\dd y\le t^{-d}\int e^{-U}$. Hence $\E_{\mu_x}e^{U/2}\le2^d$. Also $\E_\varphi e^{|v|^2/4}=2^{d/2}$. Markov's inequality applied to $e^{H/2}$ concludes.
\end{proof}

\subsection{Moments of the unadjusted chain on the region where friction dominates}\label{sec:moments}

By the change of variables $x\mapsto\sqrt m\,x$, $t\mapsto\sqrt m\,t$, which leaves $H$, the chain and the quantities $\bar h$, $\kappa$, $l_H$ and $t_E$ unchanged, we assume in Sections~\ref{sec:moments}, \ref{sec:comparison} and~\ref{sec:proofmain} that $m=1$, and also that $x_*=0$ and $U_*=0$. Then $h=\bar h/\sqrt\kappa$, $\gamma=2$, $r=e^{-h}$, $L_H=l_H$, and we assume $\bar h\le1/8$, so $h\le1/8$ and $2h(1-h)\le\sigma^2\le2h$. We use: $|x|^2\le2U(x)$, $\langle x,\nabla U(x)\rangle\ge U(x)+\frac12|x|^2$, $|\nabla U|^2\le2\kappa U$ and $G\le\frac{\bar h^2}4U\le U/256$. Put
\[
\Ly(x,v)=U(x)-G(x)+\tfrac12|v|^2+\tfrac15\langle x,v\rangle,\qquad\Emod=H-G,\qquad v_*=\frac{0.24\,\kappa}{l_H\bar h^2}=\frac{0.24}{L_Hh^2},
\]
$\delta_H=L_Hh^3/12$, so that $\delta_Hv_*=0.02h$ when $L_H>0$, and $\Dst=\{H\le v_*^2/32\}$. When $L_H=0$ we set $v_*=\infty$ and $\Dst=\R^{2d}$, and $\tau\equiv0$. The threshold $v_*$ is chosen so that the cubic Verlet error is bounded by $0.02h|V|^2$, which can be absorbed by the dissipative terms. Since $\frac15|\langle x,v\rangle|\le\frac15H$ and $G\le H/256$,
\begin{equation}\label{eq:equiv}
0.79H\le\Ly\le1.2H,\qquad0.996H\le\Emod\le H.
\end{equation}
For $w=(x,v)$ also $|v-\frac h2\nabla U(x)|\le|v|+\frac{\bar h}{\sqrt2}U^{1/2}\le1.01\sqrt{2H(w)}$ and $|V(w)|^2\le2|v|^2+\bar h^2U(x)\le4H(w)$.

\begin{lemma}[Drift]\label{lem:drift}
Let $z=(x,v)\in\Dst$, let $Z_1$ be one unadjusted step from $z$ with noises $\xi,\xi'$, $V$ its flight velocity, and $N_g=\{\sigma|\xi|\le v_*/2\}$ ($N_g=\Omega$ if $v_*=\infty$). On $N_g$, $|V|\le0.76v_*$ and $\tau\le0.02h|V|^2$. With $\theta_0=1/20$,
\[
\E\big[e^{\theta_0\Ly(Z_1)}1_{N_g}\big]\le\exp\big(\theta_0(1-0.076h)\Ly(z)+2.04\,\theta_0hd\big).
\]
\end{lemma}

\begin{proof}
Write $v'=rv+\sigma\xi$, $g_0=\nabla U(x)$, $V=v'-\frac h2g_0$, $x_1=x+hV$, $g_1=\nabla U(x_1)$, $v'''=V-\frac h2g_1$ and $Z_1=(x_1,rv'''+\sigma\xi')$. Let $a=\frac15$.

\emph{Exact expansion.} An Ornstein-Uhlenbeck half-step $(y,u)\mapsto(y,ru+\sigma\zeta)$ changes $\Ly$ by $\frac12(r^2-1)|u|^2+a(r-1)\langle y,u\rangle+\mathcal N(y,u,\zeta)$ with $\mathcal N=r\sigma\langle u,\zeta\rangle+a\sigma\langle y,\zeta\rangle+\frac{\sigma^2}2|\zeta|^2$. $\Phi$ changes $\Emod$ by $\tau$ (Lemma~\ref{lem:energy}), and $\langle x_1,v'''\rangle-\langle x,v'\rangle=h|V|^2-\frac h2[\langle x,g_0\rangle+\langle x_1,g_1\rangle]$, since $\langle x_1,v'''\rangle=\langle x,V\rangle+h|V|^2-\frac h2\langle x_1,g_1\rangle$ and $\langle x,V\rangle=\langle x,v'\rangle-\frac h2\langle x,g_0\rangle$.

\emph{Deterministic parts.} Since $|\langle y,u\rangle|\le\frac14|y|^2+|u|^2$, $1-r\le h$ and $1-r^2\ge2h(1-h)$, an Ornstein-Uhlenbeck half-step contributes at most $-h(1-h-a)|u|^2+\frac{ah}4|y|^2+\mathcal N$. With $-\frac{ah}2\langle y,\nabla U(y)\rangle\le-\frac{ah}2U(y)-\frac{ah}4|y|^2$ for $y=x,x_1$, the terms in $|x|^2$ and $|x_1|^2$ cancel, and
\[
\Ly(Z_1)-\Ly(z)\le-h(1-h-a)(|v|^2+|v'''|^2)-\tfrac{ah}2(U(x)+U(x_1))+\tau+ah|V|^2+\mathcal N_1+\mathcal N_2,
\]
with $\mathcal N_1=\mathcal N(x,v,\xi)$ and $\mathcal N_2=\mathcal N(x_1,v''',\xi')$.

\emph{The energy error.} On $\Dst$, $|rv-\frac h2g_0|\le1.01\sqrt{2H(z)}\le0.2525v_*$, and on $N_g$, $\sigma|\xi|\le v_*/2$, so $|V|\le0.76v_*$. By Lemma~\ref{lem:energy}, $\tau\le\delta_H|V|^3\le\delta_Hv_*|V|^2=0.02h|V|^2=\frac{ah}{10}|V|^2$. Since $|V|^2\le2|v'''|^2+\bar h^2U(x_1)$, $\tau+ah|V|^2\le2.2ah|v'''|^2+1.1ah\bar h^2U(x_1)$. With $h,\bar h\le\frac18$ this gives, on $N_g$,
\[
\Ly(Z_1)-\Ly(z)\le\mathcal M:=-0.675h|v|^2-0.1hU(x)-0.235h|v'''|^2-0.0965hU(x_1)+\mathcal N_1+\mathcal N_2 .
\]
\emph{Gaussian integrals.} For $\zeta\sim\mathcal N(0,\Id)$, $b\in\R^d$ and $s<1$, $\E e^{\langle b,\zeta\rangle+\frac s2|\zeta|^2}=(1-s)^{-d/2}e^{|b|^2/(2(1-s))}$. Now $\mathcal N_2=\langle b_2,\xi'\rangle+\frac{\sigma^2}2|\xi'|^2$ with $b_2=\sigma(rv'''+ax_1)$ and $|b_2|^2\le4h(|v'''|^2+0.08U(x_1))$. With $s=\theta_0\sigma^2\le1/80$,
\[
\E[e^{\theta_0\mathcal N_2}\,|\,\xi]\le\exp\big(1.02\theta_0hd+\theta_0h(0.1013|v'''|^2+0.0082U(x_1))\big),
\]
and these terms are absorbed by $-0.235h|v'''|^2-0.0965hU(x_1)$. In the same way, with $b_1=\sigma(rv+ax)$, $\E e^{\theta_0\mathcal N_1}\le\exp(1.02\theta_0hd+\theta_0h(0.1013|v|^2+0.0082U(x)))$. Hence
\begin{multline*}
\E[e^{\theta_0\Ly(Z_1)}1_{N_g}]\le e^{\theta_0\Ly(z)}\E e^{\theta_0\mathcal M}\\\le\exp\big(\theta_0\Ly(z)-\theta_0h(0.5737|v|^2+0.0918U(x))+2.04\theta_0hd\big),\end{multline*}
and $0.5737|v|^2+0.0918U\ge0.0918H\ge0.0765\Ly$ by \eqref{eq:equiv}.
\end{proof}

Let $z$ satisfy $H(z)\le R$ with $R\le v_*^2/32$, let $(Z_k)$ be the unadjusted chain from $z$, and define the events
\[
B_0=\Omega,\qquad B_{k+1}=B_k\cap\{Z_k\in\Dst\}\cap\{\sigma|\xi_k|\le v_*/2\},
\]
so that $B_k$ is determined by the first $k$ steps and, by Lemma~\ref{lem:drift}, $|V_j|\le0.76v_*$ for $j<k$ on $B_k$. Put $Q_0=1.2R+27d$ and $Q=Q_0+50$.

\begin{lemma}[Stopped moments and exits]\label{lem:stopped}
Let $z$, $(Z_k)$, $B_k$, $Q_0$ and $Q$ be as above. Then:
\begin{enumerate}
\item $\E[e^{\theta_0\Ly(Z_k)}1_{B_k}]\le e^{\theta_0Q_0}$ for all $k\ge0$.
\item For $1\le\omega\le\frac32$: $\E[1_{B_k}H(Z_k)^\omega]\le(1.266Q)^\omega$.
\item $\PP(B_n^c)\le n\,[\exp(-\theta_0(0.79v_*^2/32-Q_0))+\PP(\sigma|\xi|>v_*/2)]$, and if $v_*\ge2\sqrt d$ the second term is at most the first.
\item For $s\ge0$: $\PP(\max_{k\le n}H(Z_k)>s)\le\PP(B_n^c)+(n+1)\exp(-\theta_0(0.79s-Q_0))$. In particular, with $R''=1.27(Q_0+20\log(100(n+1)))$, $\PP(\max_{k\le n}H(Z_k)>R'')\le\PP(B_n^c)+0.01$.
\item For $k<n$: $\E[|V_k|^31_{B_{k+1}}]\le(6.6Q)^{3/2}$.
\end{enumerate}
\end{lemma}

\begin{proof}
(i) Let $y_k=\E[e^{\theta_0\Ly(Z_k)}1_{B_k}]$. By Lemma~\ref{lem:drift} applied conditionally on the first $k$ steps, $y_{k+1}\le e^{2.04\theta_0hd}\E[1_{B_k\cap\{Z_k\in\Dst\}}e^{\theta_0(1-0.076h)\Ly(Z_k)}]$. By H\"older's inequality, $\E[1_AY^{1-0.076h}]\le(\E[1_AY])^{1-0.076h}$, so
\[
\log y_{k+1}\le2.04\theta_0hd+(1-0.076h)\log y_k .
\]
Since $2.04/0.076\le27$, this recursion preserves $\log y_k\le\theta_0Q_0$, and $\log y_0=\theta_0\Ly(z)\le1.2\theta_0R$.

(ii) Let $Y=1_{B_k}e^{\theta_0\Ly(Z_k)}+1_{B_k^c}\ge1$, so $\E Y\le2e^{\theta_0Q_0}$. The function $y\mapsto(\log y)^\omega$ is concave on $[e^{\omega-1},\infty)$, so by Jensen's inequality $\E[1_{B_k}(\theta_0\Ly)^\omega]\le\E[(\log(e^\omega Y))^\omega]\le(\omega+\log2+\theta_0Q_0)^\omega$. Divide by $\theta_0^\omega$, use $20(\omega+0.7)\le50$ and \eqref{eq:equiv}.

(iii) We have
\[
B_k\setminus B_{k+1}\subset\big(B_k\cap\{H(Z_k)>v_*^2/32\}\big)\cup\{\sigma|\xi_k|>v_*/2\} .
\]
On $\{H>v_*^2/32\}$, $\Ly>0.79v_*^2/32$, and Markov's inequality with (i) bounds the first probability. Sum over $k<n$. Since $\sigma\le1/2$ and $|\xi|$ is $1$-Lipschitz with $\E|\xi|\le\sqrt d$, $\PP(\sigma|\xi|>v_*/2)\le\PP(|\xi|>v_*)\le e^{-(v_*-\sqrt d)^2/2}\le e^{-v_*^2/8}$ when $v_*\ge2\sqrt d$, and $v_*^2/8\ge\theta_0\cdot0.79v_*^2/32$.

(iv) $\{H(Z_k)>s\}\subset B_n^c\cup(B_k\cap\{\Ly(Z_k)>0.79s\})$ for $k\le n$; use (i) and Markov's inequality. For $s=R''$, $0.79\cdot1.27\ge1$ gives $\theta_0(0.79R''-Q_0)\ge\log(100(n+1))$.

(v) On $B_{k+1}$, $|V_k|\le1.01\sqrt{2H(Z_k)}+\sigma|\xi_k|$, with $\xi_k$ independent of $Z_k$. With $(a+b)^3\le4(a^3+b^3)$, (ii) and $\E|\sigma\xi|^3\le(2h(d+1))^{3/2}\le(0.01Q)^{3/2}$, $\E[|V_k|^31_{B_{k+1}}]\le4[(2.583Q)^{3/2}+(0.01Q)^{3/2}]\le(6.6Q)^{3/2}$.
\end{proof}

\subsection{A comparison chain}\label{sec:comparison}

At flight speeds above $v_*$ the energy error of velocity Verlet is controlled only by $|\tau|\le\frac{\bar h^2}8|V|^2$, which allows the energy to grow at a rate of order $\bar h^2/h=\bar h\sqrt\kappa$ per unit time. We do not control the unadjusted chain there. Instead we compare it with a chain that reverses the velocity in place of a flight faster than $v_*$. Let $S(x,v)=(x,-v)$, $F=\{(x,v):|v|>v_*\}$ and $S_F(x,v)=(x,-v)$ on $F$, $S_F=\mathrm{id}$ on $F^c$. Define
\[
\Phi_\sharp(w)=\begin{cases}\Phi(w),&|V(w)|\le v_*,\\ S(w),&|V(w)|>v_*,\end{cases}\qquad K_\sharp(z,\cdot)=\Law\big((O_{h/2}^{\xi'}\circ\Phi_\sharp\circ O_{h/2}^{\xi})(z)\big).
\]
The comparison chain is used only for estimates. It follows velocity Verlet when the proposed flight speed is at most $v_*$, and otherwise reverses the velocity without moving the position. Driven by the same noises, it agrees with the unadjusted chain until the first flight faster than $v_*$.

\begin{lemma}\label{lem:sharp}
The map $\Phi_\sharp$ has the following properties. (i) $\Phi_\sharp=B_{h/2}\circ A_{h/2}\circ S_F\circ A_{h/2}\circ B_{h/2}$; it is a Borel bijection with Borel inverse and preserves Lebesgue measure.
(ii) $H(\Phi_\sharp(w))-H(w)=G(x_1^\sharp)-G(x)+\tau_\sharp(w)$, where $x_1^\sharp$ is the position of $\Phi_\sharp(w)$ and $\tau_\sharp=\tau\,1\{|V(w)|\le v_*\}\le0.02h|V(w)|^2$.
(iii) $\Emod(\Phi_\sharp(w))\le(1+0.081h)\Emod(w)$.
\end{lemma}

\begin{proof}
(i) $B_{h/2}(w)=(x,V(w))$. If $|V|\le v_*$, then $A_{h/2}\circ S_F\circ A_{h/2}$ acts as $A_h$. If $|V|>v_*$, it maps $(x,V)$ to $(x,-V)$, and $B_{h/2}(x,-V)=(x,-V-\frac h2\nabla U(x))=S(w)$. Each factor is a Borel bijection preserving Lebesgue measure; $A_{h/2}$ preserves $F$ and $F^c$, and $S_F$ is an involution.
(ii) On $\{|V|\le v_*\}$ use Lemma~\ref{lem:energy}; on the complement $H\circ S=H$ and the position does not move. The bound is $\delta_H|V|^3\le\delta_Hv_*|V|^2$.
(iii) $|V|^2\le4H(w)\le4.02\,\Emod(w)$ and $0.02\times4.02\le0.081$.
\end{proof}

The comparison chain from $z$ is $Z_0^\sharp=z$, $Z^\sharp_1=Z_1$, and $W^\sharp_k=O_{h/2}^{\xi_k}(Z^\sharp_k)$, $Z^\sharp_{k+1}=O_{h/2}^{\xi_k'}(\Phi_\sharp(W^\sharp_k))$ for $k\ge1$, with the same noises as $(Z_k)$; $V^\sharp_k=V(W^\sharp_k)$. The first step is a step of $K$ so that the law of $Z^\sharp_k$ has a bounded density. On $B_k$ we have $|V_j|\le v_*$ for $j<k$, hence $Z^\sharp_j=Z_j$ for $j\le k$, and
\begin{equation}\label{eq:TVsharp}
\mathrm{TV}(\delta_zK^E,\Law(Z^\sharp_E))\le\PP(B_E^c).
\end{equation}

The next lemma bounds $\Ev$ along a flight in which part of the mass reverses at mid-flight. For $\alpha\in[0,1]$ let $\mathrm H_2(\alpha)=-\alpha\log\alpha-(1-\alpha)\log(1-\alpha)$.

\begin{lemma}[Conditioning on the speed]\label{lem:split}
Let $p$ be a law on $\R^{2d}$ with finite second moments and finite entropy, $\alpha=p(F)\in(0,1)$, $a=p(\cdot|F^c)$ and $b=p(\cdot|F)$. Then
\[
(1-\alpha)\Ev(a)+\alpha\Ev(b)-\mathrm H_2(\alpha)\le\Ev(p)\le(1-\alpha)\Ev(a)+\alpha\Ev(b).
\]
\end{lemma}

\begin{proof}
Let $\nu,\nu_a,\nu_b$ be the position marginals and $\varphi_F=\varphi(|v|>v_*)$. Since $F$ is determined by $v$, the chain rule for relative entropy with respect to the indicator of $F$ gives
\[
\Ev(p)=\KL\big(\mathrm{Bern}(\alpha)\,\|\,\mathrm{Bern}(\varphi_F)\big)+(1-\alpha)\KL(a\,\|\,\nu\otimes\varphi(\cdot|F^c))+\alpha\KL(b\,\|\,\nu\otimes\varphi(\cdot|F)).
\]
Here $\KL(b\|\nu\otimes\varphi(\cdot|F))=\Ev(b)+\KL(\nu_b\|\nu)+\log\varphi_F$, and similarly for $a$. The Bernoulli term plus $(1-\alpha)\log(1-\varphi_F)+\alpha\log\varphi_F$ equals $-\mathrm H_2(\alpha)$, and $(1-\alpha)\KL(\nu_a\|\nu)+\alpha\KL(\nu_b\|\nu)$ is the mutual information between the indicator of $F$ and the position, which lies in $[0,\mathrm H_2(\alpha)]$.
\end{proof}

\begin{proposition}[Corrector bound for one step of $\Phi_\sharp$]\label{prop:DSsharp}
Let $p'$ be a law with finite second moments, finite entropy and a density bounded by $\|p'\|_\infty$, let $V=V(w)$ for $w\sim p'$, $F_V=\{|V|>v_*\}$ and $\alpha=p'(F_V)\le e^{-1}$. Then \eqref{eq:DS} holds with $\Phi_\sharp$ in place of $\Phi$, $a_1=\frac{13}2$, $a_2=\frac12$ and
\begin{align*}
b_\sharp(p')={}&\tfrac{13}4\bar h^2\E|V|^2+\tfrac{13}2\E[|V|^21_{F_V}]+\tfrac2h\big(\E[|V|^21_{F_V}]\big)^{1/2}W_2(\nu_{h/2},\mu_x)\\
&+7\Big[3\alpha\log\tfrac1\alpha+2\alpha+\alpha\log^+\|p'\|_\infty+\alpha d\log(2\pi)+\E\big[(2U(x)+|V|^2)1_{F_V}\big]\Big],
\end{align*}
where $\nu_{h/2}$ is the position marginal at mid-flight.
\end{proposition}

\begin{proof}
If $\alpha=0$, then $\Phi_\sharp=\Phi$ $p'$-a.s., and the claim follows from Proposition~\ref{prop:DS}, since $b(p')\le b_\sharp(p')$. Let $\alpha>0$, so that $\alpha\in(0,1)$. Let $p_0=p'B_{h/2}$, $p_s=p_0A_s$ for $s\le h/2$, $p_{h/2}^+=p_{h/2}S_F$, $p_s=p^+_{h/2}A_{s-h/2}$ for $s\in(h/2,h]$, so that $p'\Phi_\sharp=p_hB_{h/2}$ by Lemma~\ref{lem:sharp}(i). The velocity $u$ satisfies $|u|=|V|$ throughout, so $p_s(F)=\alpha$ for all $s$. $S_F$ fixes positions and preserves $\nu\otimes\varphi$ for every $\nu$, so $\Ev(p^+_{h/2})=\Ev(p_{h/2})$, the position marginal does not jump, and $\Ev(p_s)+\Leb(\nu_s)$ is constant on $[0,h]$ by \eqref{eq:EvLeb}. Moreover $\Ccor(p^+_{h/2})-\Ccor(p_{h/2})=-2\E_{p_{h/2}}[\langle u,x-T(x)\rangle1_F]$, which is at most $\mathcal F:=2(\E[|V|^21_{F_V}])^{1/2}W_2(\nu_{h/2},\mu_x)$ in absolute value.

\emph{Corrector.} Apply Lemma~\ref{lem:kick} at both ends and Lemma~\ref{lem:approx} on $[0,h/2]$ from $p_0$ and on $[h/2,h]$ from $p_{h/2}^+$. As in Proposition~\ref{prop:DS0}, and using that $\Ev+\Leb$ is constant across the reversal, we obtain \eqref{eq:star} for $\Phi_\sharp$ with the additional term $\mathcal F$ on the right.

\emph{Free-flight bound.} Write $p_s=(1-\alpha)a_s+\alpha b_s$ with $a_s=p_s(\cdot|F^c)$ and $b_s=p_s(\cdot|F)$. The slow part is not reversed, so $a_s=a_0A_s$ and, by Corollary~\ref{cor:freeflight}, $\Ev(a_s)\le\Ev(a_0)+\Ev(a_h)$. By Lemma~\ref{lem:split} at $s$, $0$ and $h$,
\[
\Ev(p_s)\le(1-\alpha)[\Ev(a_0)+\Ev(a_h)]+\alpha\Ev(b_s)\le\Ev(p_0)+\Ev(p_h)+2\mathrm H_2(\alpha)+\alpha\Ev(b_s).
\]
Since $\KL(b_s\|\mu)=\Ev(b_s)+\KL(\Law_{b_s}(x)\|\mu_x)$, $\alpha\Ev(b_s)\le\alpha\KL(b_s\|\mu)=\int_Fp_s\log(p_s/\mu)+\alpha\log(1/\alpha)$. The density of $p_s$ is bounded by $\|p'\|_\infty$, since it is the image of $p'$ under Lebesgue-preserving bijections, $\log(1/\mu)=H+\log\int e^{-H}\le H+d\log(2\pi)$, and at the state $(x+s'V,\pm V)$ reached from $w=(x,v)$ at time $s'\in[0,h]$ of the flight, $H\le U(x+s'V)+\frac12|V|^2\le2U(x)+\bar h^2|V|^2+\frac12|V|^2\le2U(x)+|V|^2$, since $U(x+z)\le2U(x)+\kappa|z|^2$ by $L$-smoothness and $|\nabla U|^2\le2\kappa U$ (here $m=1$ and $U_*=0$). With $\mathrm H_2(\alpha)\le\alpha\log(1/\alpha)+\alpha$, this gives
\[
\int_0^h\Ev(p_s)\dd s\le h\big[\Ev(p_0)+\Ev(p_h)\big]+h\,e_{\rm mid},
\]
where $e_{\rm mid}$ is the bracket in $b_\sharp$.

\emph{Kicks.} As in Lemma~\ref{lem:identities}(ii), $\Ev(p_0)+\Ev(p_h)=\Ev(p')+\Ev(p'\Phi_\sharp)+\varrho_0-\varrho_1$. On $F_V^c$ the contribution to $\varrho_0-\varrho_1$ is at most $\frac{\bar h^2}2|V|^2$ as before. On $F_V$ the final velocity before the second kick is $-V$ at the position $x$, so the contribution is $-h\langle V,\nabla U(x)\rangle-2G(x)\le|V|^2$. Hence $\varrho_0-\varrho_1\le\frac{\bar h^2}2\E|V|^2+\E[|V|^21_{F_V}]$.

Combining as in the proof of Proposition~\ref{prop:DS} gives the claim, with $\mathcal F/h$ in $b_\sharp$.
\end{proof}

\begin{proposition}[Entropy decay of the comparison chain]\label{prop:KLsharp}
Let $\sqrt mh\le1/8$ and $n\ge1$. Then
\begin{align*}
\KL(\Law Z^\sharp_{n+1}\|\mu)\le{}&\tfrac{17}{15}e^{-\lambda nh}\KL(\delta_zK\|\mu)+\tfrac{16}{15}\Big[\tfrac{\bar h^2}4\E U(X^\sharp_{n+1})\\&+\tfrac1{\lambda h}\sup_{1\le k\le n}\E\tau_\sharp(W^\sharp_k)^++4.79\sup_{1\le k\le n}b_\sharp(\Law W^\sharp_k)\Big],
\end{align*}
provided each $\Law W^\sharp_k$ has finite second moments, finite entropy, a bounded density and $\PP(|V^\sharp_k|>v_*)\le e^{-1}$, as Proposition~\ref{prop:DSsharp} requires.
\end{proposition}

\begin{proof}
Proposition~\ref{prop:R} with $\Phi'=\Phi_\sharp$ (Lemma~\ref{lem:sharp}), $p_0=\delta_zK$ and $\theta=r/16$, as in the proof of Theorem~\ref{thm:KL}, with Proposition~\ref{prop:DSsharp} in place of Proposition~\ref{prop:DS}; $\eta/\lambda=(16+r)/(4r)\le4.79$ for $r\ge e^{-1/8}$. The laws $\Law Z^\sharp_k$ have finite moments of all orders by Lemma~\ref{lem:sharp}(iii), and finite entropy.
\end{proof}

\begin{lemma}[Moments of the comparison chain]\label{lem:momsharp}
Let $H(z)\le R$ and $t_E=Eh$. For $1\le k\le E$,
\begin{align*}
\|\Emod(Z^\sharp_k)\|_{L^2}&\le1.008e^{0.081kh}\big(R+2(d+2)kh\big),\\ \|\Emod(W^\sharp_k)\|_{L^2}&\le\Xi_E:=1.01e^{0.081t_E}\big(R+2(d+2)(t_E+1)\big).
\end{align*}
The density of $\Law W^\sharp_k$ is bounded by $(2\pi h\sigma^2)^{-d}e^{2dt_E}$.
\end{lemma}

\begin{proof}
An Ornstein-Uhlenbeck half-step from $(y,u)$ replaces $\varpi=\frac12|u|^2$ by $\varpi'=\frac12|ru+\sigma\zeta|^2$, with $\E\varpi'=r^2\varpi+\frac{\sigma^2}2d$ and $\E\varpi'^2=r^4\varpi^2+(2+d)r^2\sigma^2\varpi+\frac{\sigma^4}4(d^2+2d)$. With $A=\Emod-\varpi\ge0$ and $\sigma^2\le2h$, $\E(A+\varpi')^2\le(\Emod+h(d+2))^2$, so the $L^2$ norm of $\Emod$ increases by at most $h(d+2)$. A step of $\Phi$ multiplies $\Emod$ by at most $1+\frac{\bar h^2}8\cdot4.02\le1.008$ (Lemma~\ref{lem:energy}), and a step of $\Phi_\sharp$ by at most $\rho=1+0.081h$. So $a_k=\|\Emod(Z^\sharp_k)\|_2$ satisfies $a_1\le1.008(a_0+b)+b$ and $a_{k+1}\le\rho(a_k+b)+b$ with $b=h(d+2)$, and induction gives $a_k\le1.008\rho^{k-1}(a_0+2kb)$. Add $b$ for $W^\sharp_k$. For the density: the density of $\delta_zK$ is bounded by $(2\pi h\sigma^2)^{-d}$ (Lemma~\ref{lem:start}); Lebesgue-preserving bijections do not change the supremum, and an Ornstein-Uhlenbeck half-step multiplies it by at most $r^{-d}=e^{dh}$, since $\int\mathcal N(v';ru,\sigma^2\Id)\dd u=r^{-d}$.
\end{proof}

\section{Proof of the main theorem}\label{sec:proofmain}

We first establish mixing under four explicit conditions on the step size and the epoch length. These conditions control rejection, numerical error, entropy decay, and escape from the region where friction dominates. We then choose the parameters to satisfy all four conditions and count the resulting gradient evaluations.

\begin{proposition}[Explicit conditions]\label{prop:explicit}
Assume (A1) and (A2), and let $\varepsilon$, $\nu$, $M_2$ and $N$ be as in Theorem~\ref{thm:main}. Put
\[
R=3d\log2+2\log\big(16NM_2^2/\varepsilon^2\big),\qquad Q=1.2R+27d+50.
\]
Let $\bar h>0$ and $E\ge2$, and put $v_*=0.24\kappa/(l_H\bar h^2)$ ($v_*=\infty$ if $L_H=0$) and
\[
\Psi=1.01e^{0.081t_E}\big(R+2(d+2)(t_E+1)\big)+Q+d\log(\kappa/\bar h^2)+2dt_E .
\]
Suppose that
\begin{enumerate}
\item[(a)] $\bar h^2\big[Q+20\log(100(E+1))\big]\le1/5000$;
\item[(b)] $l_H\bar h^2(6.6Q)^{3/2}(1+t_E)\le\kappa/1000$;
\item[(c)] $\lambda(E-1)h\ge\log\big(227\,[1.02R+d\log(\kappa/\bar h^2)]\big)$;
\item[(d)] $v_*^2\ge41\big[Q+20\log(2E)+80\log(\Psi\sqrt\kappa/\bar h)+922\big]$.
\end{enumerate}
Then $\mathrm{TV}(\nu P^N,\mu)\le\varepsilon$.
\end{proposition}

Condition (a) controls the bias term of the entropy decay (Proposition~\ref{prop:KLsharp}) and the term $\frac{\bar h^2}4\max_nU(X_n)$ in the rejection probability of an epoch; neither contains $\kappa$. Condition (b) controls the trapezoidal rule errors in both places. Condition (c) makes the epoch long enough for the relative entropy to decay. Condition (d) makes the probability that the unadjusted chain leaves the region where friction dominates the energy error, within one epoch, smaller than any power of the parameters that appears in the proof; since (b) gives $v_*\ge240(6.6Q)^{3/2}(1+t_E)$, condition (d) holds unless $\Psi\sqrt\kappa/\bar h$ is exponentially large in $Q^3(1+t_E)^2$. With fresh uniforms, the rejection term would instead be of order $E\,\bar h^3Q$.

In the proof we use the conventions of Section~\ref{sec:moments}: $m=1$, $x_*=0$, $U_*=0$. Condition (a) gives $\bar h^2\le1/(5000\cdot50)$, so $\bar h\le1/8$ and $h\le1/8$.

\begin{proof}[Proof of Proposition~\ref{prop:explicit}]
\emph{Reduction.} Let $\mathcal D=\{H\le R\}$. By Lemma~\ref{lem:tail} and the choice of $R$, $N\mu(\mathcal D^c)\le\varepsilon^2/(16M_2^2)$, and $M_2\ge1$. By Proposition~\ref{prop:coupling}, $\mathrm{TV}(\nu P^N,\mu)\le q^N+\varepsilon/16+\varepsilon/4$, and $q^N\le\varepsilon/2$ if $q\le1/4$. By the triangle inequality through $\mu$, it suffices to show that for every $z\in\mathcal D$
\[
\mathrm{TV}(\delta_zK^E,\mu)\le0.105\qquad\text{and}\qquad\beta(z)\le0.0125 .
\]
Indeed, these give $q\le2(0.105+0.0125)=0.235<\frac14$, and then $\mathrm{TV}(\nu P^N,\mu)\le\frac\varepsilon2+\frac\varepsilon{16}+\frac\varepsilon4<\varepsilon$.

Fix $z\in\mathcal D$ and put $p_*=\PP(B_E^c)$, with the events $B_k$ of Section~\ref{sec:moments}.

\emph{The exit probability.} This bounds the cost of replacing the unadjusted chain by the comparison chain. Condition (d) gives $v_*\ge2\sqrt d$ and $0.79v_*^2/32\ge Q_0+20\log(2E)+80\log(\Psi/h)+921.1$, so Lemma~\ref{lem:stopped}(iii) gives
\[
p_*\le2E\exp\big(-\theta_0(0.79v_*^2/32-Q_0)\big)\le10^{-20}(h/\Psi)^4 .
\]
It also gives $R\le v_*^2/32$, so $\mathcal D\subset\Dst$ and Section~\ref{sec:moments} applies. Note $\Psi\ge\Xi_E+Q+d\log(\kappa/\bar h^2)+2dt_E\ge1$, so $\Xi_Ep_*^{1/2}\le10^{-10}h^2$.

\emph{Rejection.} This establishes the second target, $\beta(z)\le0.0125$. Let $R''=1.27(Q_0+20\log(100(E+1)))$ and $\mathcal G=B_E\cap\{\max_{n\le E}H(Z_n)\le R''\}$. By Proposition~\ref{prop:allaccept}, $U\le H$ and Lemma~\ref{lem:stopped}(iv),
\[
\beta(z)\le\PP(\mathcal G^c)+\tfrac{\bar h^2}4R''+\sum_{k<E}\E\big[|\tau(W_k)|1_{B_{k+1}}\big]\le2p_*+0.01+\tfrac{\bar h^2}4R''+E\delta_H(6.6Q)^{3/2},
\]
using $B_E\subset B_{k+1}$, $|\tau(W_k)|\le\delta_H|V_k|^3$ and Lemma~\ref{lem:stopped}(v). By (a), $\frac{\bar h^2}4R''\le\frac{1.27}{4\cdot5000}$, and $E\delta_H=t_El_H\bar h^2/(12\kappa)$, so by (b) the last term is at most $1/12000$. Hence $\beta(z)\le0.0102$.

\emph{Moments of the comparison chain.} These feed the bias terms of Proposition~\ref{prop:KLsharp}. For $1\le k\le E$, by Lemma~\ref{lem:stopped}(ii), Lemma~\ref{lem:momsharp}, \eqref{eq:equiv} and $Z^\sharp_k=Z_k$ on $B_k$,
\[
\E H(Z^\sharp_k)\le\E[H(Z_k)1_{B_k}]+1.005\,\Xi_Ep_*^{1/2}\le1.27Q .
\]
For $1\le k<E$, since $\xi_k$ is independent of $Z^\sharp_k$, $\E|V^\sharp_k|^2\le2.0402\,\E H(Z^\sharp_k)+2hd\le2.61Q$, and $\||V^\sharp_k|^2\|_{L^2}\le4.02\Xi_E$. The event $F_V=\{|V^\sharp_k|>v_*\}$ is contained in $B_{k+1}^c$, so $\alpha:=\PP(F_V)\le p_*$. On $B_{k+1}$, $V^\sharp_k=V_k$, so by Lemma~\ref{lem:sharp}(ii) and Lemma~\ref{lem:stopped}(v),
\[
\E\tau_\sharp(W^\sharp_k)^+\le\delta_H(6.6Q)^{3/2}+0.02h\,\E[|V^\sharp_k|^21_{B_{k+1}^c}]\le\delta_H(6.6Q)^{3/2}+10^{-11}h .
\]

\emph{The remainder $b_\sharp$.} This bounds the additional terms produced by the velocity reversals of the comparison chain. We bound the terms of Proposition~\ref{prop:DSsharp} at $p'=\Law W^\sharp_k$. First $\E[|V|^21_{F_V}]\le4.02\Xi_E\alpha^{1/2}\le4.1\cdot10^{-10}h^2$ and $\E[(2U+|V|^2)1_{F_V}]\le6.1\Xi_E\alpha^{1/2}$. Next, $\E_{\mu_x}|x|^2\le d$ (integrate $\langle x,\nabla U\rangle\ge|x|^2$ by parts), and the mid-flight position $x+\frac h2V$ satisfies $\E|x+\frac h2V|^2\le4\E U(X^\sharp_k)+\frac{h^2}2\E|V|^2\le5.1Q$, so $W_2(\nu_{h/2},\mu_x)\le2.26\sqrt Q+\sqrt d\le2.5\sqrt Q$. Hence the third term of $b_\sharp$ is at most $\frac{5\sqrt Q}h(4.02\Xi_E)^{1/2}p_*^{1/4}\le10.1\cdot10^{-5}(Q\Xi_E)^{1/2}/\Psi\le5.1\cdot10^{-5}$. By Lemma~\ref{lem:momsharp}, $\log^+\|p'\|_\infty+d\log(2\pi)\le d\log(\kappa/\bar h^2)+2dt_E\le\Psi$, using $h\sigma^2\ge1.75\bar h^2/\kappa$. With $\alpha\log(1/\alpha)\le\alpha^{1/2}$, the bracket in $b_\sharp$ is at most $p_*^{1/2}(5+\Psi p_*^{1/2}+6.1\Xi_E)\le10^{-9}$, and the second term is at most $3\cdot10^{-9}$. So $b_\sharp\le\frac{13}4\bar h^2\E|V^\sharp_k|^2+6\cdot10^{-5}$, and $\alpha\le e^{-1}$.

\emph{Relative entropy.} This establishes the first target. By Lemma~\ref{lem:start} and $\E_zH(Z_1)\le1.013(R+hd)+hd\le1.02R+0.27d$, $\KL(\delta_zK\|\mu)\le1.02R+d\log(\kappa/\bar h^2)$, using $\log(1/(2e(1-\sqrt mh)))\le-1.56$, valid since $\sqrt mh\le \bar h\le0.002$ by (a). Proposition~\ref{prop:KLsharp} with $n=E-1$, $(\lambda h)^{-1}\le86.8/h$ and the bounds above give
\begin{align*}
\KL(\Law Z^\sharp_E\|\mu)\le{}&\tfrac{17}{15}e^{-\lambda(E-1)h}\big(1.02R+d\log(\kappa/\bar h^2)\big)\\&+\tfrac{16}{15}\Big[40.95\,\bar h^2Q+7.234\,\frac{l_H\bar h^2}{\kappa}(6.6Q)^{3/2}+3\cdot10^{-4}\Big],
\end{align*}
since $\frac{1.27}4+4.79\cdot\frac{13}4\cdot2.61\le40.95$ and $86.8\,\delta_H/h\le7.234\,l_H\bar h^2/\kappa$. By (c) the first term is at most $0.005$; by (a) and (b) the second is at most $\frac{16}{15}(0.00819+0.00724+0.0003)\le0.0168$. So
\[
\KL(\Law Z^\sharp_E\|\mu)\le0.005+0.0168=0.0218,
\]
and by Pinsker's inequality and \eqref{eq:TVsharp},
\[
\mathrm{TV}(\delta_zK^E,\mu)\le p_*+(0.0218/2)^{1/2}\le0.105 .
\]
Hence $q\le0.21+0.0204\le1/4$.
\end{proof}

\begin{proof}[Proof of Theorem~\ref{thm:main}]
Let $R$ and $Q$ be as in Proposition~\ref{prop:explicit}. Since $M_2\ge1$ and $\varepsilon<1$, $D\ge d\ge1$, and $\log N\le\log(2/\varepsilon)\le1+\log(1/\varepsilon)$, so $R\le6D+8$ and $D\le Q\le94D$.

\emph{Choice of parameters.} We choose an epoch duration $t_0=C\Lambda$, long enough for the entropy term to decay. The step size must satisfy two restrictions: one controls the leading energy error, the other the accumulated trapezoidal remainder over the epoch. We therefore take the smaller of the two corresponding step sizes, and then choose the number of steps so that the epoch lasts at least $t_0$. Let $C\ge\mathrm e$ be a constant to be fixed, put $t_0=C\Lambda$, and for $\bar h\in(0,1/8]$ let $E=\lceil t_0\sqrt\kappa/\bar h\rceil+1$, so that $t_0\le t_E\le t_0+1$. Take
\[
\bar h^2=\min\Big\{\frac{1}{C^2Q\Lambda},\ \frac{\kappa}{C^2l_HQ^{3/2}t_0}\Big\}.
\]
Since $Q\le94D\le94e^\Lambda$, we have $\log(\kappa/\bar h^2)\le C_1\Lambda\log C$, $\log E\le C_1\Lambda\log C$ and $\log\Psi\le0.081(t_0+1)+C_1\Lambda\log C$ for an absolute constant $C_1$. Below, $C_1,\dots,C_5$ are absolute constants independent of $C$.

\emph{Condition (a).} Since $\bar h^2\le1/(C^2Q\Lambda)$,
\[
\bar h^2\big[Q+20\log(100(E+1))\big]\le\frac{Q+C_2\Lambda\log C}{C^2Q\Lambda}\le\frac1{5000}
\]
once $C^2\ge5000(2+C_2\log C)$.

\emph{Condition (b).} Since $\bar h^2\le\kappa/(C^2l_HQ^{3/2}t_0)$ and $1+t_E\le2+t_0$,
\[
l_H\bar h^2(6.6Q)^{3/2}(1+t_E)\le\frac{17\kappa(2+t_0)}{C^2t_0}\le\frac{\kappa}{1000}
\]
once $C^2\ge51000$.

\emph{Condition (c).} We have $\lambda\ge1/86.8$ and $(E-1)h\ge t_0$, while
\[
\log\big(227[1.02R+d\log(\kappa/\bar h^2)]\big)\le C_3\Lambda\log C,
\]
so (c) holds once $C\ge86.8C_3\log C$.

\emph{Condition (d).} We have $v_*^2\ge(0.24C^2Q^{3/2}t_0)^2\ge0.05C^6Q^3\Lambda^2$, while the right side of (d) is at most
\[
41\big[Q+922+C_4\Lambda(C+\log C)\big],
\]
so (d) holds once $C^5\ge C_5$. A single sufficiently large absolute constant $C$ satisfies all four conditions simultaneously, and Proposition~\ref{prop:explicit} gives $\mathrm{TV}(\nu P^N,\mu)\le\varepsilon$.

\emph{The gradient count.} With this $C$,
\begin{align*}
E\le\frac{t_0\sqrt\kappa}{\bar h}+2&\le C\Lambda\big(C\sqrt{\kappa Q\Lambda}+C\sqrt{l_HQ^{3/2}C\Lambda}\big)+2\\&\le C_0'\Lambda^{3/2}\big(\sqrt{\kappa D}+\sqrt{l_H}\,D^{3/4}\big)
\end{align*}
by $Q\le94D$. Also $N\le1+\log(2/\varepsilon)/\log4\le3\log(2/\varepsilon)/\log4$, since $\log4\le2\log(2/\varepsilon)$. So the bound on $NE$ holds with $C_0=3C_0'/\log4$, and the cost is $NE+1$ gradient evaluations.
\end{proof}

\section{Bayesian logistic regression}\label{sec:numerics}

\emph{The posterior.} Let $(x_i,y_i)$, $i=1,\dots,n$, be observations with covariates $x_i\in\R^d$, one entry of which is $1$ (an intercept), and labels $y_i\in\{0,1\}$. With the logistic likelihood and the prior $\mathcal N(0,s^2\Id)$, the posterior density of the coefficients $\theta\in\R^d$ is
\begin{equation}\label{eq:posterior}
\pi(\theta)=\frac1Z\exp\Big(-\sum_{i=1}^n\Big[\log\big(1+e^{\langle x_i,\theta\rangle}\big)-y_i\langle x_i,\theta\rangle\Big]-\frac{|\theta|^2}{2s^2}\Big),
\end{equation}
where $Z$ is the normalizing constant. Neither $Z$ nor posterior expectations, such as the posterior mean of $\theta$ or the posterior predictive probability $\int(1+e^{-\langle x,\theta\rangle})^{-1}\pi(\theta)\,\dd\theta$ of the label $1$ at covariates $x$, are available in closed form. They are computed by sampling from $\pi$, and logistic regression is a standard benchmark for this task \cite{CR17}.

\emph{The assumptions.} Let $\psi(t)=\log(1+e^t)$, so that $0<\psi''\le\frac14$ and $|\psi'''|\le\frac1{6\sqrt3}$. In the coordinates $\phi=\theta/s$, up to an additive constant, $\pi$ has potential
\[
\begin{aligned}
U(\phi)&=\tfrac12|\phi|^2+\sum_{i=1}^n\Big[\psi\big(s\langle x_i,\phi\rangle\big)-s\,y_i\langle x_i,\phi\rangle\Big],\\
\nabla^2U(\phi)&=\Id+s^2\sum_{i=1}^n\psi''\big(s\langle x_i,\phi\rangle\big)\,x_ix_i^\top .
\end{aligned}
\]
Let $X$ be the $n\times d$ matrix with rows $x_i^\top$. Then (A1) holds with
\[
m=1,\qquad L=1+\tfrac{s^2}4\lambda_{\max}(X^\top X),
\]
and since $|\psi''(a)-\psi''(b)|\le|a-b|/(6\sqrt3)$, (A2) holds with
\[
L_H=\frac{s^3}{6\sqrt3}\,\max_i|x_i|\,\lambda_{\max}(X^\top X).
\]
The change of variables $\theta\mapsto\theta/s$ is the rescaling $x\mapsto\sqrt m\,x$ of Section~\ref{sec:moments} with $m=s^{-2}$, which leaves the chain unchanged. So Theorem~\ref{thm:main} applies to $\pi$, with $\kappa=L$, and the prior scale $s$ sets the condition number. As $s$ grows, $\kappa$ grows like $s^2$ and $L_H$ like $s^3$, so the second term of Theorem~\ref{thm:main}, $\sqrt{l_H}\,D^{3/4}$, grows like $\kappa^{3/4}$ and dominates over our range of $s$. On these posteriors Theorem~\ref{thm:main} therefore guarantees cost of order $\kappa^{3/4}$; the square-root growth measured below goes beyond what it proves.

\emph{Data.} We use two standard data sets. The breast cancer data \cite{WMSS93} have $n=569$ observations and $30$ covariates, which we standardize to mean $0$ and standard deviation $1$; with the intercept, $d=31$. The Sonar data \cite{SG88} have $n=208$ observations and $60$ covariates, which we rescale to mean $0$ and standard deviation $\frac12$ as in \cite{CR17}; with the intercept, $d=61$. In both, the data are nearly linearly separable along some directions of parameter space; there the likelihood is almost flat, the curvature of $U$ is set by the prior, and the condition number grows with $s$.

\emph{Chains.} We compare adjusted OBABO with a fresh uniform at every step and with the non-reversibly updated uniform (Definition~\ref{alg:pu}), redrawn once per unit of time, with $\gamma=2$. We run $32$ chains of $1000$ units of time, started near the mode after a burn-in of $60$ units, and report the integrated autocorrelation time (IAT) in units of time, the maximum over the coordinates, $|\phi|^2$ and $U$, estimated from block means with a self-consistent window. Each step costs one gradient evaluation, so the number of gradient evaluations per effective sample is the IAT divided by $h$.

\emph{At a step based on the global curvature.} We first take $h=\bar h/\sqrt L$ with the practical value $\bar h=2/\sqrt{d+2}$ instead of the constants of Theorem~\ref{thm:main}, and $16$ chains of $400$ units of time. For $s=0.2,0.8,3.2$ (breast cancer) and $s=0.25,1,4$ (Sonar), that is $\kappa$ from $11$ to $1.9\cdot10^4$, both chains reject fewer than $0.2\%$ of the moves, and all IATs lie between $3.1$ and $4.6$ units of time. The IAT in units of time does not grow with $\kappa$; at these rejection rates the two acceptance rules give nearly the same chain.

\emph{At the curvature of the mode.} The global bound $L$ exceeds the largest eigenvalue $\hat L$ of $\nabla^2U$ at the mode by a factor $8$ to $52$ (breast cancer) and $1.4$ to $12$ (Sonar) over our range of $s$, so the step size of Theorem~\ref{thm:main} is conservative. We therefore also run both chains with $h=1/\sqrt{\hat L}$, which is outside the parameter range of Theorem~\ref{thm:main}. Since the smallest eigenvalue of $\nabla^2U$ at the mode lies between $1$ and $1.35$, $\hat\kappa=\hat L$ is, up to this factor, the condition number of the posterior at its mode. Both chains then reject approximately $7\%$ to $11\%$ (breast cancer) and $7\%$ to $13\%$ (Sonar) of the moves. We vary $s$ from $0.2$ to $51.2$ (breast cancer) and from $0.25$ to $64$ (Sonar), so that $\hat\kappa=\hat L$ ranges from $9.5$ to $9.5\cdot10^4$ and from $8$ to $5.6\cdot10^4$. For Sonar at $s=128$ this step gives a rejection rate of $25\%$, outside the range above, and we do not include it. The curvature varies along the chain: the largest eigenvalue of $\nabla^2U$ lies between $0.38\hat L$ and $1.22\hat L$ at $s=12.8$ (breast cancer) and between $0.59\hat L$ and $0.89\hat L$ at $s=16$ (Sonar), as $10\%$ and $90\%$ quantiles. So the chains do not see a nearly Gaussian target, for which \cite{AGS24} proves a $\sqrt\kappa$ bound.

\begin{figure}[t]
\centering
\includegraphics[width=\textwidth]{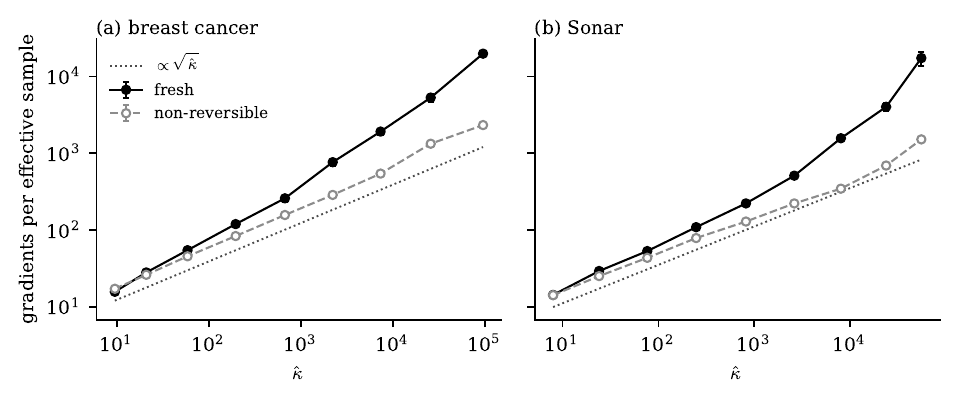}
\caption{Gradient evaluations per effective sample against $\hat\kappa$ for the (a) breast cancer and (b) Sonar posteriors.}
\label{fig:hard}
\end{figure}

Figure~\ref{fig:hard} shows the cost. Over these ranges, the IAT of the chain with the non-reversibly updated uniform stays approximately between $5.6$ and $8.3$ (breast cancer) and between $3.9$ and $6.4$ (Sonar) units of time, while that of the fresh-uniform chain grows from $5.1$ to $64$ and from $5.1$ to $74$. Least-squares lines fitted to the logarithm of the cost against $\log\hat\kappa$ have slopes $0.54\pm0.01$ (non-reversibly updated uniform) and $0.76\pm0.02$ (fresh uniform) for breast cancer, and $0.50\pm0.02$ and $0.76\pm0.05$ for Sonar, with standard errors from the residuals. The scripts \texttt{logistic\_global.py} and \texttt{logistic\_sweep.py}, in the ancillary files of the arXiv version of this paper, with \texttt{LOGIT\_DATA} set to \texttt{breast\_cancer} or \texttt{sonar}, reproduce these runs.

\section{Discussion and open problems}\label{sec:discussion}

Unless stated otherwise, the chain considered here is adjusted OBABO with the non-reversibly updated uniform of Definition~\ref{alg:pu}.

\emph{The Hessian-Lipschitz condition.} In our argument, $L_H$ enters through the trapezoidal error $\tau$ of Lemma~\ref{lem:energy}. A Lipschitz Hessian gives the cubic bound $|\tau|\le\frac{L_H}{12}h^3|V|^3$. Under strong convexity and a Lipschitz gradient alone, the uniform bound is quadratic in $h|V|$. This order is sharp: for $d=1$, the potential equal to $\frac{a+\Delta}2x^2$ for $x\le0$ and to $\frac a2x^2$ for $x\ge0$ gives $|\tau|=\frac\Delta8h^2V^2$ on a flight from $-\frac h2V$ to $\frac h2V$, and smooth approximations preserve the curvature bounds and approach this value. A H\"older-continuous Hessian of exponent $\alpha$ instead gives an error of order $h^{2+\alpha}|V|^{2+\alpha}$. Whether the chain itself needs a Lipschitz Hessian for square-root acceleration remains open. A sharper rejection bound could exploit cancellation between signed errors, but must control their maximum partial sum over an epoch. The experiments also motivate sharper bounds: the stated $L_H$ bound for the logistic posteriors makes the second term of Theorem~\ref{thm:main} grow like $\kappa^{3/4}$, whereas the measured cost grows approximately like $\hat\kappa^{1/2}$. Here $\hat\kappa$ is the condition number at the mode, and the experiments use a different step size and measure the stationary IAT rather than total variation.

\emph{Cold starts.} Theorem~\ref{thm:main} applies whenever the initial density ratio has finite $L^2(\mu)$ norm $M_2$. When $\log M_2=O(d\log(2\kappa))$, its first complexity term is $\tilde O(\sqrt{\kappa d})$. For a poorly initialized density with $\log M_2$ of order $\kappa d$, this term becomes $\tilde O(\kappa\sqrt d)$. A point mass has no density ratio, so the theorem does not apply to it directly. Remark~\ref{rem:warm} supplies an initialization procedure: from $x_0$, accelerated gradient descent finds an approximate minimizer using $O\big(\sqrt\kappa\log(2+\kappa L|x_0-x_*|^2)\big)$ gradient evaluations, after which a Gaussian draw gives $M_2\le e(2\kappa)^{d/2}$. These initialization gradients are additional to the sampling cost. An accelerated bound for the adjusted chain started directly from an arbitrary point, without this optimization stage, remains open.

\emph{Dependence on the dimension.} On Gaussian targets, the first term of Theorem~\ref{thm:main} gives $\tilde O(\sqrt\kappa\,d^{1/2})$, whereas Proposition~\ref{prop:sharp} gives a lower bound proportional to $\sqrt\kappa\,d^{1/6}$. The upper bound controls the position-dependent energy term $G=\frac{h^2}8|\nabla U|^2$ uniformly over an energy sublevel set whose level is at least of order $d$. Improving the dimension dependence requires sharper control of rejection and larger step sizes. Whether this chain attains $\tilde O(\sqrt\kappa\,d^{1/4})$, as obtained in \cite{AGS24} for Gaussian targets and in \cite{CCLZ26,YGC26} for the full class under their warm-start assumptions, is open. The entropy estimate also leaves room for improvement: on Gaussian targets its moment remainder is of order $\bar h^2d$ at stationarity, whereas the stationary relative entropy of unadjusted OBABO on Gaussian targets is $O(\bar h^4d)$ for steps bounded away from the stability boundary $\bar h=2$.

\emph{Other integrators.} Theorem~\ref{thm:KL} is proved for OBABO. Whether an analogous entropy estimate yields accelerated mixing for other Metropolis-adjusted splittings, including OABAO with the non-reversibly updated uniform, remains open.

\appendix

\section{Proof of the free-transport estimate \NoCaseChange{\eqref{eq:freesmooth}}}\label{app:free}

In this appendix $\xi$ denotes the displacement $y-T(y)$ of a Brenier map, as in \cite{Lu26}, and not a Gaussian noise. Throughout, $\tilde p$, $\iota$, $q$, $q_s$ and $\nu_s$ are as in \eqref{eq:freesmooth}, $t>0$ is fixed, $|z|\le R_0$ on the support of $\tilde p$, and $C$ denotes constants that do not depend on $s\in[0,t]$ or $y\in\R^d$; they may change from line to line. For a law $q'$ on $\R^{2d}$ and $(Y,W)\sim q'$ write $m(y)=\E[W|Y=y]$, $\Sigma(y)=\mathrm{Cov}(W|Y=y)$ and $\pi_y$ for the conditional law of $W$ given $Y=y$, so that $J(q')=\E|m(Y)|^2$ and, by the chain rule, $\Ev(q')=\E\KL(\pi_Y\|\varphi)$. For a matrix field $S$, $(\mathrm{div}\,S)_i=\sum_j\partial_jS_{ij}$, $(\nabla m)_{ij}=\partial_jm_i$ and $A:B=\tr(AB^\top)$.

\subsection{The velocity covariance term}

\begin{lemma}\label{lem:stress}
Assume that $\mu_x$ has finite second moments and $\E_{\mu_x}|U|<\infty$. Let $q'$ be a law on $\R^{2d}$ with finite second moments and finite entropy whose position marginal $\nu$ has a continuous, strictly positive density $n$, let $T$ be the Brenier map from $\nu$ to $\mu_x$ and $\xi(y)=y-T(y)$. Assume that $n\Sigma$ is $C^1$ and that, for some $C$, $|n\Sigma|\le Cn$ and $|\mathrm{div}(n\Sigma)|\le C(1+|y|)n$. Then
\[
-\int\langle\xi,\mathrm{div}(n\Sigma)\rangle\,\dd y\le4\Ev(q')-2J(q')-\Delta\Leb(\nu).
\]
\end{lemma}

In the notation of \cite{Lu26}, the left side is $\mathfrak S-\mathfrak A+\int q\nabla U\cdot\xi\,\dd\mu_x$, and, without the term $-2J$, the lemma is what the proof of \cite[Lemma~5.3]{Lu26} gives together with the identity in the proof of Lemma~\ref{lem:kick}, before the convexity of $U$ is used there. We do not split $\Sigma$ into $\Id$ and $\Sigma-\Id$; then the singular part of the Hessian of the Brenier potential enters once, with a sign, whereas Lu treats the two parts separately. The term $-2J$ comes from centering the conditional laws $\pi_y$.

\begin{proof}
\emph{The Brenier potential.} By Brenier's theorem \cite{Bre91}, $T=\nabla\phi$ $\nu$-a.e.\ for a convex $\phi:\R^d\to\R\cup\{+\infty\}$. The convex set $\{\phi<\infty\}$ has full $\nu$-measure, hence full Lebesgue measure because $n>0$, so its interior is $\R^d$, and $\phi$ is finite and locally Lipschitz. The distributional Hessian of $\phi$ is a nonnegative definite symmetric matrix of Radon measures, $D^2\phi=G\,\dd y+D^s\phi$ with $D^s\phi$ singular to Lebesgue measure, and $G$ is the Alexandrov Hessian of $\phi$, which exists a.e.\ \cite{EG15}. Since $\nu$ and $\mu_x$ have densities, McCann's theorem \cite{McC97} gives, for a.e.\ $y$, $\det G(y)>0$ and $n(y)=Z_x^{-1}e^{-U(T(y))}\det G(y)$ with $Z_x=\int e^{-U}$. As $\int n|\log n|<\infty$ and $\int n|U\circ T|\,\dd y=\int|U|\,\dd\mu_x<\infty$, $\log\det G$ is $\nu$-integrable and
\begin{equation}\label{eq:MA}
\int n\log\det G\,\dd y=\Leb(\nu)+\E_{\mu_x}U+\log Z_x=\Leb(\nu)-\Leb(\mu_x)=\Delta\Leb(\nu).
\end{equation}

\emph{Integration by parts.} Put $S=n\Sigma$, which is nonnegative definite. Let $\chi\in C_c^\infty(\R^d)$ be radial and nonincreasing in $|y|$, with $\chi=1$ on the unit ball and $0\le\chi\le1$, and put $\chi_k(y)=\chi(y/k)$. For a $C^1$ symmetric matrix field $\Theta$ with compact support, $\int\langle\nabla\phi,\mathrm{div}\,\Theta\rangle\,\dd y=-\int\Theta:\dd D^2\phi$ and $\int\langle y,\mathrm{div}\,\Theta\rangle\,\dd y=-\int\tr\Theta\,\dd y$. With $\Theta=\chi_kS$ and $\mathrm{div}(\chi_kS)=\chi_k\,\mathrm{div}\,S+S\nabla\chi_k$ this gives
\[
-\int\chi_k\langle\xi,\mathrm{div}\,S\rangle\,\dd y=\int\langle\xi,S\nabla\chi_k\rangle\,\dd y+\int\chi_kS:(\Id-G)\,\dd y-\int\chi_kS:\dd D^s\phi .
\]
The last integral is nonnegative, because $S$ and $D^s\phi$ are nonnegative definite. Since $\int n|\xi|^2=W_2^2(\nu,\mu_x)<\infty$ and $\nu$ has finite second moments, the assumptions on $S$ and dominated convergence give $\int\chi_k\langle\xi,\mathrm{div}\,S\rangle\to\int\langle\xi,\mathrm{div}\,S\rangle$, $|\int\langle\xi,S\nabla\chi_k\rangle|\le Ck^{-1}\int n|\xi|\to0$ and $\int\chi_k\tr S\to\int\tr S<\infty$, while $\int\chi_kS:G\uparrow\int S:G$ by monotone convergence, since $S:G\ge0$. Hence $\int S:G<\infty$ and
\begin{equation}\label{eq:ibp}
-\int\langle\xi,\mathrm{div}\,S\rangle\,\dd y\le\int n\,\Sigma:(\Id-G)\,\dd y .
\end{equation}

\emph{Gaussian duality.} Fix $y$ with $\det G(y)>0$ and put $K=\Id-G(y)\preceq\Id$. Let $\pi_y^c$ be $\pi_y$ translated to mean zero; then $\KL(\pi^c_y\|\varphi)=\KL(\pi_y\|\varphi)-\frac12|m(y)|^2$. The Gibbs variational principle with the function $w\mapsto\frac14\langle w,Kw\rangle$, whose exponential has $\varphi$-integral $\det(\Id-K/2)^{-1/2}$, gives
\[
\tfrac14\,\Sigma(y):K\le\KL(\pi_y\|\varphi)-\tfrac12|m(y)|^2-\tfrac12\log\det(\Id-K/2).
\]
If $g_1,\dots,g_d>0$ are the eigenvalues of $G(y)$, then $-2\log\det(\Id-K/2)=-2\sum_i\log\frac{1+g_i}2\le-\sum_i\log g_i$, since $(1+g)^2\ge4g$. Therefore $\Sigma(y):(\Id-G(y))\le4\KL(\pi_y\|\varphi)-2|m(y)|^2-\log\det G(y)$. Integrate against $n$ and use \eqref{eq:ibp} and \eqref{eq:MA}.
\end{proof}

\subsection{Free transport of $q$}

Let $(X,V)\sim\tilde p$ and $\zeta=(\zeta_x,\zeta_v)\sim\mathcal N(0,\Id_{2d})$ be independent. For $s\in\R$, $q_s$ is the law of $(Y_s,W)$ with
\[
W=V+\sqrt\iota\,\zeta_v,\qquad Y_s=X+sV+\sqrt\iota(\zeta_x+s\zeta_v).
\]
Put $\sigma_s^2=\iota(1+s^2)$ and $\beta_s=s/(1+s^2)$. Given $(X,V)$, $Y_s\sim\mathcal N(X+sV,\sigma_s^2\Id)$, and given $(X,V,Y_s)$, $W\sim\mathcal N\big(V+\beta_s(Y_s-X-sV),\iota(1+s^2)^{-1}\Id\big)$. We write $n_s$, $m_s$, $\Sigma_s$ for the density of $\nu_s$ and the conditional mean and covariance of $W$ given $Y_s$.

\emph{Densities.} The density of $q_s$ is the Gaussian convolution $z'\mapsto\E\,g_s(z'-A_s(X,V))$, where $g_s$ is the density of $\mathcal N(0,\iota A_sA_s^\top)$, so it is $C^\infty$ in $(s,z')$ and each of its derivatives is bounded by a polynomial in $|z'|$ times $e^{-|z'|^2/C}$, uniformly for $s\in[0,t]$. In the same way $n_s(y)=\E\,\mathcal N(y;X+sV,\sigma_s^2\Id)$ is $C^\infty$, and since $|X+sV|\le C$,
\begin{equation}\label{eq:nbounds}
C^{-1}e^{-C(1+|y|^2)}\le n_s(y)\le Ce^{-|y|^2/C}\qquad(s\in[0,t],\ y\in\R^d).
\end{equation}

\emph{Moment equations.} Since $q_s(y,w)=q(y-sw,w)$, $\partial_sq_s=-\langle w,\nabla_yq_s\rangle$. With $J_s=\int wq_s(\cdot,w)\,\dd w=n_sm_s$ and $M_s=\int w\otimes w\,q_s(\cdot,w)\,\dd w=n_s(\Sigma_s+m_s\otimes m_s)$, differentiation under the integral gives
\begin{equation}\label{eq:moments}
\partial_sn_s=-\mathrm{div}\,J_s,\qquad\partial_sJ_s=-\mathrm{div}\,M_s .
\end{equation}

\emph{The conditional mean velocity.} The conditional law of $(X,V)$ given $Y_s=y$ is $\pi_{s,y}(\dd z)\propto e^{\psi(s,y,z)}\tilde p(\dd z)$, $z=(x,u)$, with $\psi(s,y,z)=\langle y,x+su\rangle/\sigma_s^2-|x+su|^2/(2\sigma_s^2)$. Write $\langle\cdot\rangle_{s,y}$ for expectation under $\pi_{s,y}$ and $F_s(x,u)=(1+s^2)^{-1}u-\beta_sx$. Since $\E[W|X,V,Y_s]=\beta_sY_s+F_s(X,V)$,
\[
m_s(y)=\beta_sy+b(s,y),\quad b(s,y)=\langle F_s\rangle_{s,y},\qquad\Sigma_s(y)=\frac{\iota}{1+s^2}\Id+\mathrm{Cov}_{\pi_{s,y}}(F_s),
\]
so $|\Sigma_s|\le C$. We claim that $|\partial_s^j\partial_y^\alpha b(s,y)|\le C_{j,\alpha}(1+|y|)^j$ on $[0,t]\times\R^d$ for all $j$ and multi-indices $\alpha$. Indeed, for $F$ smooth in $s$ and bounded with its $s$-derivatives on $[0,t]\times\{|z|\le R_0\}$, and $\partial\in\{\partial_s,\partial_{y_i}\}$, $\partial\langle F\rangle_{s,y}=\langle\partial F\rangle_{s,y}+\langle F\,\partial\psi\rangle_{s,y}-\langle F\rangle_{s,y}\langle\partial\psi\rangle_{s,y}$. By induction, $\partial_s^j\partial_y^\alpha\langle F\rangle_{s,y}$ is a finite sum of products of averages $\langle F_1\rangle_{s,y}\cdots\langle F_l\rangle_{s,y}$, where each $F_i$ is a product of derivatives of $F$ and of $\psi$, and at most $j$ of these factors are of the form $\partial_s^k\psi$, $k\ge1$. Each such factor is affine in $y$ with bounded coefficients, $\partial_y\partial_s^k\psi$ is bounded and independent of $y$, and $\partial_y^2\psi=0$. This proves the claim. Consequently, on $[0,t]\times\R^d$,
\begin{gather*}
|m|\le C(1+|y|),\qquad|\nabla m|+|\nabla^2m|\le C,\\
|\partial_sm|+|\nabla\partial_sm|\le C(1+|y|),\qquad|\partial_s^2m|\le C(1+|y|)^2 .
\end{gather*}

\emph{The acceleration.} Put $a_s=\partial_sm_s+(\nabla m_s)m_s$ and $e_s=\partial_sa_s+(\nabla a_s)m_s$. The bounds above give $|a_s|\le C(1+|y|)$ and $|e_s|\le C(1+|y|)^2$. By \eqref{eq:moments}, $n_s\partial_sm_s=\partial_sJ_s-m_s\partial_sn_s=-\mathrm{div}\,M_s+m_s\,\mathrm{div}\,J_s$, and $\mathrm{div}(J_s\otimes J_s/n_s)=m_s\,\mathrm{div}\,J_s+n_s(\nabla m_s)m_s$, so
\[
n_sa_s=-\mathrm{div}\big(M_s-J_s\otimes J_s/n_s\big)=-\mathrm{div}(n_s\Sigma_s).
\]
Thus $n_s\Sigma_s$ is $C^\infty$, $|n_s\Sigma_s|\le Cn_s$ and $|\mathrm{div}(n_s\Sigma_s)|\le C(1+|y|)n_s$. Let $T_s$ be the Brenier map from $\nu_s$ to $\mu_x$ and $\xi_s(y)=y-T_s(y)$. Lemma~\ref{lem:stress} applies to $q_s$ for every $s\in[0,t]$ and gives
\begin{equation}\label{eq:Qbound}
\int\langle a_s,\xi_s\rangle\,n_s\,\dd y\le4\Ev(q_s)-2J(q_s)-\Delta\Leb(\nu_s).
\end{equation}

\subsection{Coupling along the flow of the mean velocity}

Since $m$ is smooth and $|\nabla m|\le C$, the equation $\dot y=m_{s'}(y)$ has a global flow $(\Phi_{s,s'})_{s,s'\in[0,t]}$ of diffeomorphisms. By \eqref{eq:moments}, $(\nu_s)$ is a narrowly continuous solution of $\partial_s\nu_s+\mathrm{div}(m_s\nu_s)=0$, and since $m$ is Lipschitz in $y$, $\nu_{s'}=(\Phi_{s,s'})_\#\nu_s$ \cite[Proposition~8.1.8]{AGS08}, for $s'\ge s$ and, by the group property $\Phi_{s',s}=\Phi_{s,s'}^{-1}$, also for $s'<s$; the integrability conditions there hold by $|m|\le C(1+|y|)$, $|\nabla m|\le C$ and \eqref{eq:nbounds}. Let $\mathcal W(s)=\frac12W_2^2(\nu_s,\mu_x)$. Fix $s\in[0,t]$ and $\vartheta$ with $s+\vartheta\in[0,t]$, and put $X_\vartheta=\Phi_{s,s+\vartheta}$. Then $\partial_\vartheta X_\vartheta=m_{s+\vartheta}(X_\vartheta)$, $\partial_\vartheta^2X_\vartheta=a_{s+\vartheta}(X_\vartheta)$ and $\partial_\vartheta^3X_\vartheta=e_{s+\vartheta}(X_\vartheta)$, and $|X_\vartheta(y)|\le C(1+|y|)$ by Gronwall's lemma. By Taylor's formula,
\[
X_\vartheta(y)=y+\vartheta m_s(y)+\tfrac12\vartheta^2a_s(y)+\rho_\vartheta(y),\qquad|\rho_\vartheta(y)|\le C|\vartheta|^3(1+|y|)^2 .
\]
For $Y\sim\nu_s$ the pair $(X_\vartheta(Y),T_s(Y))$ is a coupling of $\nu_{s+\vartheta}$ and $\mu_x$, so $\mathcal W(s+\vartheta)\le\frac12\E|\xi_s(Y)+X_\vartheta(Y)-Y|^2$. Expand the square. Since $\E|\xi_s(Y)|^2=2\mathcal W(s)\le C$ and, by \eqref{eq:nbounds}, the moments of $\nu_s$ are bounded uniformly in $s$, the Cauchy-Schwarz inequality gives
\begin{equation}\label{eq:Wexp}
\mathcal W(s+\vartheta)\le\mathcal W(s)+\vartheta\,\Ccor(q_s)+\tfrac12\vartheta^2Q(s)+M|\vartheta|^3,\qquad Q(s)=J(q_s)+\int\langle a_s,\xi_s\rangle\,n_s\,\dd y,
\end{equation}
with $M$ independent of $s$ and $\vartheta$. Here we used $\E\langle m_s(Y),\xi_s(Y)\rangle=\E\langle W,\xi_s(Y_s)\rangle=\Ccor(q_s)$ and $\E|m_s(Y)|^2=J(q_s)$.

\subsection{Conclusion}

Let $N\ge1$, $\vartheta=t/N$ and $s_i=i\vartheta$. Apply \eqref{eq:Wexp} at $s_i$ with increment $\vartheta$ and at $s_{i+1}$ with increment $-\vartheta$, and add:
\[
\Ccor(q_{s_{i+1}})-\Ccor(q_{s_i})\le\tfrac\vartheta2\big[Q(s_i)+Q(s_{i+1})\big]+2M\vartheta^2 .
\]
By \eqref{eq:Qbound}, $Q\le\bar Q$ with $\bar Q(s)=4\Ev(q_s)-J(q_s)-\Delta\Leb(\nu_s)$. By \eqref{eq:EvLeb}, $4\Ev(q_s)-\Delta\Leb(\nu_s)$ equals $5\,\mathrm h(\nu_s)$ plus a constant. The functions $s\mapsto\mathrm h(\nu_s)=-\int n_s\log n_s$ and $s\mapsto J(q_s)=\int|m_s|^2n_s$ are continuous on $[0,t]$ by dominated convergence, since $n_s$ and $m_s$ are continuous in $s$ and, by \eqref{eq:nbounds}, $|n_s\log n_s|+|m_s|^2n_s\le C(1+|y|^2)e^{-|y|^2/C}$. So $\bar Q$ is continuous. Summing over $i<N$,
\[
\Ccor(q_t)-\Ccor(q)\le\sum_{i<N}\tfrac\vartheta2\big[\bar Q(s_i)+\bar Q(s_{i+1})\big]+2Mt\vartheta\ \longrightarrow\ \int_0^t\bar Q(s)\,\dd s\qquad(N\to\infty).
\]
Since $-J\le J\le2\Ev$ by (S1), $\bar Q(s)\le6\Ev(q_s)-\Delta\Leb(\nu_s)$, and this proves \eqref{eq:freesmooth}. Keeping $-J$ would improve the constants of Proposition~\ref{prop:DS} slightly; we do not use this.

\section{Proof of Proposition~\ref{prop:sharp}}\label{app:sharp}

\emph{Reduction to one step at stationarity.} Let $\bar\mu=\mu$ for the fresh uniform and $\bar\mu=\mu\otimes\Exp(1)$ on states $(z,\ell)$ for the non-reversibly updated uniform. Write $A=\{x_1\le0\}$, also for the set of states whose position lies in it, so that $\bar\mu(A)=\frac12$. Let $\nu_k$ be the law of the state at the start of step $k+1$, after any refresh of $\ell$. Initially
\[
\nu_0=2\,\mathbf 1_A\,\bar\mu\le2\bar\mu .
\]
Every step, and every refresh of $\ell$, is a Markov kernel that preserves $\bar\mu$ (Lemma~\ref{lem:invariance} for the second rule), so the bound propagates: $\nu_k\le2\bar\mu$ for every $k$. This is a bound on the density of $\nu_k$ with respect to $\bar\mu$, not a claim that $\nu_k=\bar\mu$; within an epoch, $\ell$ is not $\Exp(1)$ independent of $z$ under $\nu_k$. Let $c(y)$ be the probability that one step from the state $y$ moves the position from $A$ to $A^c$; in particular $c=0$ outside $A$. The position changes only on accepted moves, so a crossing is an accepted move with $x_1\le0<x_1'$, where $x_1'$ is the first coordinate of the proposed position $x'$. Then
\[
\PP(X_k\in A,\ X_{k+1}\notin A)=\int c\,\dd\nu_k\le2\int c\,\dd\bar\mu .
\]
A step started from $\bar\mu$ itself applies $O_{h/2}$, which leaves the state distributed as $\bar\mu$; for the fresh uniform the move from $w$ is then accepted with probability $\min(1,e^{-\Delta H(w)})$, and for the non-reversibly updated uniform with probability $\PP(\ell>\Delta H(w))=\min(1,e^{-\Delta H(w)})$, since under $\bar\mu$ the variable $\ell\sim\Exp(1)$ is independent of $w$. So for both rules $\int c\,\dd\bar\mu=\bar\mu(A)\,\Phi=\frac12\Phi$, where, with $w\sim\mu$,
\[
\Phi:=\frac{1}{\mu(A)}\,\E\big[\mathbf 1_{\{x_1\le0<x_1'\}}\min(1,e^{-\Delta H(w)})\big]
\]
is the probability, at stationarity and conditionally on starting in $A$, that one step leaves $A$. Thus each crossing probability of the actual chain is at most $\Phi$. Since $X_0\in A$, the event $X_n\notin A$ requires a crossing at some step $k<n$, so
\[
\PP(X_n\notin A)\le\sum_{k=0}^{n-1}\PP(X_k\in A,\ X_{k+1}\notin A)\le n\Phi,
\]
and, taking $B=A^c$ in $\mathrm{TV}(p,q)=\sup_B|p(B)-q(B)|$,
\[
\mathrm{TV}(\Law(X_n),\mu_x)\ge\mu_x(A^c)-\PP(X_n\notin A)\ge\tfrac12-n\Phi .
\]
This is the conductance lower bound for mixing times, applied to the extended chain.

\emph{Factorization.} For quadratic $U$, $\tau\equiv0$ and $\Delta H(w)=G(x')-G(x)$ by Lemma~\ref{lem:energy}, a sum $\sum_{i=1}^d\Delta H_i$ of terms depending on $(x_i,v_i)$ only, and these pairs are independent under $\mu$. Using $\min(1,e^{-y})\le e^{-y/2}$,
\[
\Phi\le2\,\E\big[\mathbf 1_{\{x_1\le0<x_1'\}}e^{-\Delta H_1/2}\big]\prod_{i=2}^d\E\big[e^{-\Delta H_i/2}\big].
\]
In the coordinate $i$, with curvature $\lambda_i$, put $\omega=h\sqrt{\lambda_i}$, $q=\sqrt{\lambda_i}\,x_i$ and $p=v_i$, which are independent standard Gaussians. Then $q'=\sqrt{\lambda_i}\,x_i'=(1-\frac{\omega^2}2)q+\omega p$ and $\Delta H_i=\frac{\omega^2}8(q'^2-q^2)$. The quadratic form $\frac12(q^2+p^2)+\frac{\omega^2}{16}(q'^2-q^2)$ is positive definite, and a Gaussian integral gives the normalizing constant
\[
\rho(\omega):=\E\big[e^{-\Delta H_i/2}\big]=\Big(1+\frac{\omega^6}{64}\Big)^{-1/2}.
\]
Under the probability measure with density $e^{-\Delta H_i/2}/\rho(\omega)$, the pair $(q,q')$ is centred Gaussian with equal variances and correlation
\[
\chi(\omega)=\frac{1-\omega^2/2}{1+\omega^4/8},\qquad 1-\chi(\omega)=\frac{\omega^2(4+\omega^2)}{8+\omega^4}\le\omega^2 .
\]
The stiff coordinates have $\omega=\bar h$. For the soft coordinate, $\omega=\bar h/\sqrt\kappa$, and the quadrant formula $\PP(Y_1\le0<Y_2)=\frac1{2\pi}\arccos\chi$ for a centred Gaussian pair with correlation $\chi$ gives
\[
2\,\E\big[\mathbf 1_{\{x_1\le0<x_1'\}}e^{-\Delta H_1/2}\big]=\frac{\rho(\omega)}{\pi}\arccos\chi(\omega)\le\sqrt{\frac{1-\chi(\omega)}2}\le\frac{\bar h}{\sqrt{2\kappa}},
\]
using $\rho\le1$ and $\arccos\chi=2\arcsin\sqrt{(1-\chi)/2}\le\pi\sqrt{(1-\chi)/2}$. Hence
\[
\Phi\le\frac{\bar h}{\sqrt{2\kappa}}\Big(1+\frac{\bar h^6}{64}\Big)^{-(d-1)/2}.
\]

\emph{Conclusion.} If $n\Phi<\frac14$, then $\mathrm{TV}(\Law(X_n),\mu_x)>\frac14$; this gives the first condition. The function $\bar h\mapsto \bar h^{-1}(1+\bar h^6/64)^{(d-1)/2}$ on $(0,\infty)$ is minimized at $\bar h^6=64/(3d-4)$, its only critical point, since it tends to $\infty$ as $\bar h\to0$ and as $\bar h\to\infty$. Velocity Verlet is stable on this Gaussian at the minimizer, since $64/(3d-4)\le32<2^6$ gives $\bar h<2$. At the minimizer the function equals $\frac12(3d-4)^{1/6}\big(\frac{3d-3}{3d-4}\big)^{(d-1)/2}\ge\frac12d^{1/6}$ for $d\ge2$. Hence the first condition holds whenever $n<\frac{\sqrt2}8\sqrt\kappa\,d^{1/6}$.

\subsection*{Acknowledgements}

I thank Bob Carpenter for introducing me to Neal's non-reversibly updated accept/reject uniform, and Pierre Monmarch\'e for discussions on acceleration from cold and warm starts. The ideas in this paper are my own. I used publicly available large language models to assist with calculations, proofs of technical lemmas, numerical experiments, figures, drafting, and the Lean formalization.

\end{document}